\documentclass[reqno,11pt]{amsart}

\newtheorem{theorem}{Theorem}[section]
\newtheorem{lemma}[theorem]{Lemma}
\newtheorem{corollary}[theorem]{Corollary}

\theoremstyle{definition}

\newtheorem{assumption}{Assumption}[section]

 \theoremstyle{remark}
\newtheorem{remark}[theorem]{Remark}

\newcommand\bB{\mathbb{B}}

\newcommand\bH{\mathbb{H}}
\newcommand\bL{\mathbb{L}}

\newcommand\bQ{\mathbb{Q}}
\newcommand\bR{\mathbb{R}}

\newcommand\bW{\mathbb{W}}

\newcommand\cB{\mathcal{B}}

\newcommand\cF{\mathcal{F}}
\newcommand\cH{\mathcal{H}}
 
\newcommand\cP{\mathcal{P}}

\newcommand\cW{\mathcal{W}}

\newcommand{\tr}{\text{\rm tr}\,}

\newcommand{\mysection}[1]{\section{#1}
\setcounter{equation}{0}}

\newcommand\cbrk{\text{$]$\kern-.15em$]$}} 
\newcommand\opar{
\text{\,\raise.2ex\hbox{${\scriptstyle |}$}\kern-.34em$($}} 
\newcommand\cpar{%
\text{$)$\kern-.34em\raise.2ex\hbox{${\scriptstyle |}$}}\,}
\newcommand\obrk{\text{$[$\kern-.15em$[$}}

\begin{document}

\title[SPDEs with VMO coefficients]
{On divergence
form SPDEs
with VMO coefficients}
\author{N.V. Krylov}
\address{127 Vincent Hall, University of Minnesota, Minneapolis,
 MN, 55455}
\thanks{The work   was partially supported by
NSF Grant DMS-0653121}
\email{krylov@math.umn.edu}

\renewcommand{\subjclassname}{%
\textup{2000} Mathematics Subject Classification}


\begin{abstract}
We present several results on
solvability in Sobolev spaces $W^{1}_{p}$
of SPDEs in divergence
form in the whole space.

\end{abstract}

\maketitle

\mysection{Introduction}

The theory of (usual) partial differential equations
has two rather different parts depending on whether
the equations are written in  divergence  or
nondivergence form.
Quite often the starting point is the same: equations with constant
coefficients, and then one uses different techniques to
treat different types of equations.

By now, one can say that the $L_{p}$-theory of evolutional 
second order SPDEs
is quite well developed.
 The most advanced results of this theory
can be found in the following papers and references therein:
\cite{Ki}
  (nondivergence
type equations), \cite{Ki1} and \cite{Ki2}
(divergence type equations). The results of the present paper
are close to the corresponding results of \cite{Ki1}. However,
unlike \cite{Ki1} we do not assume that the leading coefficients
are continuous in the space variable. Instead we assume 
that the  leading coefficients of the 
``deterministic'' part of the equation
 are in VMO which is a  much wider class
than $C$. Still the leading coefficients of the ``stochastic''
part are assumed to be continuous in~$x$.

The exposition in \cite{Ki1} and \cite{Ki2}
is based on the theory of solvability in spaces 
$H^{\gamma}_{p}=(1-\Delta)^{-\gamma/2}L_{p}$
of SPDEs with coefficients independent of $x$.
Then the method of ``freezing'' the coefficients is applied
as in the general framework set out in \cite{Kr99}.
This method does not work if the coefficients are 
only in VMO and we use a different technique based on recent
results from \cite{Kr07} on deterministic parabolic equations
with VMO coefficients. In addition, our technique allows us
to avoid using the $W^{n}_{2}$-theory
of SPDEs, which is 
a starting point in the paper \cite{Kr99} and 
subsequent articles based on it.

One more difference of our approach from the one in \cite{Ki1}
is that we represent the free term in the deterministic
part in the form $D_{i}f^{i}+f^{0}$ with $f^{j}\in L_{p}$
(see \eqref{11.13.1} below). 
Of course, this 
is just a general form of a distribution from $H^{-1}_{p}$.
However,
the spaces $H^{\gamma}_{p}$ are most appropriate 
for equations in nondivergence form.
One general inconvenience of these spaces is
that the space or space-time
dilations affect the norms in a way which
is hard to control.
For divergence form equations with low regularity of coefficients
the most important space is $H^{1}_{p}$. This space
 coincides with the Sobolev space $W^{1}_{p}$
and the effect of dilations on the norm or on
$D_{i}f^{i}+f^{0}$ can be easily
taken into account.

The exposition here is self-contained
apart from references to some very basic
results of \cite{Kr99}, \cite{Kr07}, and \cite{Trento} and
  is much more elementary than in \cite{Ki1}, employing
  the derivatives instead of the powers of the Laplacian,
 and yet gives more 
  information. In particular, the author
intends to use 
Corollary \ref{corollary 12.11.1}  in order to
largely simplify the theory in \cite{Ki1}
of divergence form SPDEs in domains.
It turns out that to develop this theory
one need not first develop the theory
of SPDEs   in domains with coefficient independent of $x$,
which in itself required quite a bit of work.

The author's interest in divergence type equations and
in simplifying the theory of them appeared after
he realized that the corresponding results
can be applied to filtering theory 
of partially observable
diffusion processes, given by stochastic It\^o equations,
and proving that,
under Lipschitz and nondegeneracy conditions only,
the filtering density is almost Lipschitz in $x$
and almost H\"older $1/2$ in time.
This is proved in \cite{Kr08_1} on the basis
of Theorems \ref{theorem 12.7.1}
through \ref{theorem 1.4.3}   of the present article.
The filtering density satisfies an SPDE usually written
in terms of the operators adjoint to operators in nondivergence
form with Lipschitz continuous coefficients.
Writing these adjoint operators in divergence form
makes perfect sense and allows us to obtain the above mentioned results
(see \cite{Kr08_1}).

Our Theorem \ref{theorem 12.7.1}
is very close to Theorem 2.12 of \cite{Ki1}.
Apart from weaker conditions on the coefficients, another
important difference is the presence of the parameter
$\lambda$ in \eqref{1.4.2}. One of differences in the proofs
is that we avoid proving  the solvability on small
consecutive time intervals
and then gluing together the results.

Let $(\Omega,\cF,P)$ be a complete probability space
with an increasing filtration $\{\cF_{t},t\geq0\}$
of complete with respect to $(\cF,P)$ $\sigma$-fields
$\cF_{t}\subset\cF$. Denote by $\cP$ the predictable
$\sigma$-field in $\Omega\times(0,\infty)$
associated with $\{\cF_{t}\}$. Let
 $w^{k}_{t}$, $k=1,2,...$, be independent one-dimensional
Wiener processes with respect to $\{\cF_{t}\}$.  

We fix a stopping time $\tau$ and for $t\leq\tau$
in the Euclidean $d$-dimensional space $\bR^{d}$
of points $x=(x^{1},...,x^{d})$ we consider the following
equation
\begin{equation}
                                          \label{11.13.1}
du_{t}=(L_{t}u_{t}-\lambda u_{t}
+D_{i}f^{i}_{t}+f^{0}_{t})\,dt
+(\Lambda^{k}_{t}u_{t}+g^{k}_{t})\,dw^{k}_{t},
\end{equation}
  where $u_{t}=u_{t}(x)=u_{t}(\omega,x)$ is an unknown function,
$$
L_{t}\psi(x)=D_{j}\big(a^{ij}_{t}(x)D_{i}\psi(x)
+a^{j}_{t}(x)\psi(x)
\big)+b^{i}_{t}(x)D_{i}\psi(x)+c_{t}(x)\psi(x),
$$
$$
\Lambda^{k}_{t}\psi(x)=\sigma^{ik}_{t}(x)D_{i}\psi(x)
+\nu^{k}_{t}(x)\psi(x),
$$
the summation  convention
with respect to $i,j=1,...,d$ and $k=1,2,...$ is enforced
and detailed assumptions on the coefficients and the 
free terms will be given later.

One can rewrite \eqref{11.13.1} in the nondivergence form
assuming that the coefficients $a^{ij}_{t}$ and $a^{j}_{t}$
are differentiable
 in $x$ and then one could apply the results
from \cite{Kr99}. It turns out that the 
differentiability of $a^{ij}_{t}$ and $a^{j}_{t}$ is not needed
for the corresponding 
counterparts of the results in \cite{Kr99} to be true
and showing this and generalizing the 
corresponding results of \cite{Ki1}
is one of the main purposes of the present
article.

The author is sincerely grateful to Kyeong-Hun Kim
who kindly pointed out an error in the first draft
of the article.

 \mysection{Main results}

Fix a number 
$$
p\geq2, 
$$
and denote $L_{p}=L_{p}(\bR^{d})$.
We use the same notation $L_{p}$ for vector- and matrix-valued
or else
$\ell_{2}$-valued functions such as
$g_{t}=(g^{k}_{t})$ in \eqref{11.13.1}. For instance,
if $u(x)=(u^{1}(x),u^{2}(x),...)$ is 
an $\ell_{2}$-valued measurable function on $\bR^{d}$, then
$$
\|u\|^{p}_{L_{p}}=\int_{\bR^{d}}|u(x)|_{\ell_{2}}^{p}
\,dx
=\int_{\bR^{d}}\big(
\sum_{k=1}^{\infty}|u^{k}(x)|^{2}\big)^{p/2}
\,dx.
$$

Introduce
$$
D_{i}=\frac{\partial}{\partial x^{i}},\quad i=1,...,d,
\quad\Delta=D^{2}_{1}+...+D^{2}_{d}.
$$
By $Du$   we mean the gradient  with respect
to $x$ of a function $u$ on $\bR^{d}$.

As usual,  
$$
W^{1}_{p}=\{u\in L_{p}: Du\in L_{p}\},
\quad
 \|u\|_{W^{1}_{p}}=
\|u\|_{L_{p}}+\|Du\|_{L_{p}}.
$$

Recall that $\tau$ is a stopping time and introduce
$$
\bL _{p}(\tau):=L_{p}(\opar 0,\tau\cbrk,\cP,
L_{p}),\quad
\bW^{1}_{p}(\tau):=L_{p}(\opar 0,\tau\cbrk,\cP,
W^{1}_{p}).
$$
We also need the space $\cW^{1}_{p}(\tau)$,
which is the space of functions $u_{t}
=u_{t}(\omega,\cdot)$ on $\{(\omega,t):
0\leq t\leq\tau,t<\infty\}$ with values
in the space of generalized functions on $\bR^{d}$
and having the following properties:

(i) We have $u_{0}\in L_{p}(\Omega,\cF_{0},L_{p})$;

(ii)  We have $u
\in \bW^{1}_{p}(\tau )$;

(iii) There exist   $f^{i}\in \bL_{p}(\tau)$,
$i=0,...,d$, and $g=(g^{1},g^{2},...)\in \bL_{p}(\tau)$
such that
 for any $\varphi\in C^{\infty}_{0}=C^{\infty}_{0}(\bR^{d})$ 
with probability 1
for all   $t\in[0,\infty)$
we have
$$
(u_{t\wedge\tau},\varphi)=(u_{0},\varphi)
+\sum_{k=1}^{\infty}\int_{0}^{t}I_{s\leq\tau}
(g^{k}_{s},\varphi)\,dw^{k}_{s}
$$
\begin{equation}
                                                 \label{1.2.1}
+
\int_{0}^{t}I_{s\leq\tau}\big((f^{0}_{s},\varphi)-(f^{i}_{s},D_{i}\varphi)
 \big)\,ds.
\end{equation}
In particular, for any $\phi\in C^{\infty}_{0}$, the process
$(u_{t\wedge\tau},\phi)$ is $\cF_{t}$-adapted and (a.s.) continuous.

The reader can find 
in \cite{Kr99} a discussion of (ii) and (iii),
in particular, the fact that the series in \eqref{1.2.1}
converges uniformly in probability on every finite
subinterval of $[0,\tau]$. On the other hand,
it is worth saying that
the above introduced space $\cW^{1}_{p}$
is not quite the same as $\cH^{1}_{p}(\tau)$ in \cite{Kr99}
or in \cite{Ki1}. 
There are 
three differences. One is that there is an additional restriction
on $u_{0}$ in \cite{Kr99} and \cite{Ki1}. But in the main part of
the article we are going to work with
$\cW^{1}_{p,0}(\tau)$ which is the subset of
$\cW^{1}_{p}(\tau)$ consisting of functions with
$u_{0}=0$. Another issue is that in \cite{Kr99} and \cite{Ki1}
we have $f^{i}=0$, $i=1,...,d$, 
and 
$$
f^{0}\in\bH^{-1}_{p}(\tau)=L_{p}(\opar0,\tau\cbrk,\cP,H^{-1}_{p}).
$$
 Actually, this 
difference is fictitious  because one knows that
any $f\in H^{-1}_{p} $ 

(a) has the form $D_{i}f^{i}+f^{0}$
with $f^{j}\in L_{p}$ and
$$
\|f\|_{H^{-1}_{p}}\leq N\sum_{j=0}^{d}\|f^{j}\|_{L_{p}},
$$
where $N$ is independent of $f,f^{j}$, and on the other hand,

(b) for any $f\in H^{-1}_{p}$  there exist
$f^{j}\in L_{p}$ such that $f=D_{i}f^{i}+f^{0}$ and
$$
\sum_{j=0}^{d}\|f^{j}\|_{L_{p}}\leq N\|f\|_{H^{-1}_{p}},
$$
where $N$ is independent of $f$.

The third difference is that instead of (i) the condition
$D^{2}u\in\bH^{-1}_{p}(\tau)$ is required in \cite{Kr99}
and \cite{Ki1}. However,
as it follows from Theorem 3.7 of \cite{Kr99}
and the boundedness of the operator $D:L_{p}\to H^{-1}_{p}$,
this difference disappears if $\tau$ is a bounded stopping time.

To summarize,
 the spaces $\cW^{1}_{p,0}(\tau)$ introduced above 
coincide with $\cH^{1}_{p,0}(\tau)$ from
 \cite{Kr99}  if $\tau$ is bounded
and we choose a particular 
representation
of the deterministic part of the stochastic differential
just for convenience.  In the remainder of the article the spaces
 $\cH^{1}_{p,0}(\tau)$ do not appear and none of their properties
is used.

In case that property (iii) holds, we write
\begin{equation}
                                       \label{12.3.1}
du_{t}=(D_{i}f^{i}_{t}+f^{0}_{t})\,dt
+g^{k}_{t}\,dw^{k}_{t}
\end{equation}
for $t\leq\tau$
and this explains the sense in which equation
\eqref{11.13.1} is understood. Of course, we still need to
specify appropriate assumptions on the coefficients
and the free terms in \eqref{11.13.1}.

\begin{assumption}
                                        \label{assumption 1.2.1}

(i) The coefficients $a^{ij}_{t}$, $a^{i}_{t}$, $b^{i}_{t}$,
$\sigma^{ik}_{t}$, $c_{t}$, and $\nu^{k}_{t}$ are measurable 
with respect to $\cP\times \cB(\bR^{d})$, where $\cB(\bR^{d})$
is the Borel $\sigma$-field on $\bR^{d}$.

(ii) There is a constant $K$ such that
for all values of indices and arguments
$$
 |a^{i}_{t}|+|b^{i}_{t}|+|c_{t}|+|\nu|_{\ell_{2}}\leq
K,\quad c_{t}\leq0.
$$

(iii) There is a constant $\delta>0$ such that
for all values of the arguments and $\xi\in\bR^{d}$
\begin{equation}
                                             \label{1.3.2}
 a^{ij}_{t}  
  \xi^{i}
\xi^{j}\leq\delta^{-1}|\xi|^{2},\quad
(a^{ij}_{t}  
-  \alpha^{ij}_{t}) \xi^{i}
\xi^{j}\geq\delta|\xi|^{2},
\end{equation}
where   $\alpha^{ij}_{t}=(1/2)(\sigma^{i\cdot},\sigma^{j\cdot})
_{\ell_{2}}$. Finally, the constant $\lambda\geq0$.

\end{assumption} 

It is worth emphasizing that we do not require
the matrix $(a^{ij})$ to be symmetric.

Assumption \ref{assumption 1.2.1}  guarantees
that equation \eqref{11.13.1} makes perfect sense
if $u\in\cW^{1}_{p}(\tau)$.
By the way, adding the term $-\lambda u_{t}$
with constant $\lambda\geq0$ is one more technically
convenient step. One can always introduce this term,
if originally it is absent, by considering $v_{t}:=u_{t}
e^{ \lambda t}$.

Let $\bB$ denote the set of balls $B\subset\bR^{d}$
and let $\rho(B)$ be the radius of $B\in\bB$.
For   functions $h_{t}(x)$ on $[0,\infty)\times\bR^{d }$
and  $B\in\bB$   introduce
$$
 h_{ t(B )} 
=\frac{1}{|B |}\int_{B }h_{t}( x)\,dx,
$$
where   $|B|$ is the volume of $B$. Also let $\bQ$ denote
the set of all cylinders in $[0,\infty)\times\bR^{d }$ of type
$Q=(s,t)\times B$, where   $B\in\bB$
and $t-s=\rho^{2}(B)$. For such $Q$ set $\rho(Q)=\rho(B)$. 
For $\rho\geq0$, $s<t$, 
a continuous $\bR^{d}$-valued function $x_{r},r\in[s,t]$,
 and a $Q=(s,t)\times B\in\bQ$, introduce
$$
\text{osc} \, (h,Q,x_{\cdot})=\frac{1}{t-s}
\int_{s}^{t}(|h_{r} -h_{ r(B+x_{r})}|)_{(B+x_{r})} \,dr,
$$
$$
\text{Osc} \, (h,Q,\rho)=\sup_{|x_{\cdot}|_{C}\leq\rho}
\text{osc} \, (h,Q,x_{\cdot}),\quad
\text{osc} \, (h,Q)=\text{osc} \, (h,Q,0),
$$
where $|x_{\cdot}|_{C}$ is the sup norm of $|x_{\cdot}|$. 
 
 Observe that $\text{osc\,}(h,Q,x_{\cdot})  =0 $ if $h_{t}(x)$ is
independent of $x$.

Denote by $B_{\rho}$ the open ball with radius $\rho>0$ centered
at the origin, define $Q_{\rho}=(0,\rho^{2})
\times B_{\rho}$ and for $t\geq0$ and $x\in\bR^{d}$
set $B_{\rho}(x)=B_{\rho}+x$, $Q_{\rho}(t,x)=Q_{\rho}+(t,x)$.

In the remaining two assumptions we use
constants $\beta>0$ and $\beta_{1}>0$ the values of which
will be specified later.

Let $t_{0}\geq0$, $x_{0}\in\bR^{d}$, and 
constants 
$\varepsilon\geq
\varepsilon_{1}>0$. 
We say that the couple $(a,\sigma)$ is $(\varepsilon,
\varepsilon_{1})$-regular
at point $(t_{0},x_{0})$ if    (for any  $\omega$)
either

(i) we have $\sigma_{t}^{nm}(x_{0})=0$ for $ t\in(t_{0},t_{0}
+\varepsilon_{1}^{2})$ 
and all $n,m$ and
\begin{equation}
                                                   \label{8.7.1}
\text{\rm osc\,}(a^{ij},Q)\leq \beta,
\quad\forall i,j,
\end{equation}
for all $Q\in\bQ$ such that $Q
 \subset Q_{\varepsilon}(t_{0},x_{0})$, or

(ii) for all $Q\in\bQ$ such that $Q
 \subset Q_{\varepsilon}(t_{0},x_{0})$  we have
\begin{equation}
                                                   \label{8.6.1}
\text{\rm Osc\,}(a^{ij},Q,\varepsilon)\leq \beta,
\quad\forall i,j.
\end{equation}

Note that $(a,\sigma)$ is $(\varepsilon,
\varepsilon_{1})$-regular
at any point $(t_{0},x_{0})$ for any $\beta>0$ if, for instance, $a^{ij}$
depend only on $x$ and are of class VMO. 
   
\begin{assumption}
                                       \label{assumption 1.2.6}
  There exist $\varepsilon\geq\varepsilon_{1}>0$ such that
$(a,\sigma)$ is $(\varepsilon,\varepsilon_{1})$-regular
at any point $(t_{0},x_{0})$ and
$$
(a^{jk}_{t}(x)  
-  \alpha^{jk}_{t}(y)) \xi^{j}\xi^{k}\geq\delta|\xi|^{2}
$$
for all   $t $, $\xi$, $x$, and $y$  satisfying
  $|x-y|\leq\varepsilon $.
\end{assumption}

\begin{assumption}
                                       \label{assumption 8.9.1} 
 There exists an $\varepsilon_{2}>0$ such that
\begin{equation}
                                                   \label{8.8.1}
 |\sigma^{i\cdot}_{t}(x)-
\sigma^{i\cdot}_{t}(y)|_{\ell_{2}}  \leq \beta_{1}
\end{equation}
for all $i$, $t $,  $x$, and $y$   satisfying
$|x-y|\leq\varepsilon_{2}$.

\end{assumption}

Needless to say that Assumptions
\ref{assumption 1.2.6} and  
\ref{assumption 8.9.1} are satisfied with
any $\beta,\beta_{1}>0$ and
slightly reduced $\delta$ if 
\eqref{1.3.2} holds 
and $a^{ij}_{t}(x)$
and $\sigma^{i\cdot}_{t}(x)$ are uniformly
continuous in $x$ uniformly with
respect to $(\omega,t)$.

Finally, we describe the space of initial data.
Recall that for $p\geq2$
the Slobodetskii space $W^{1-2/p}_{p}
=W^{1-2/p}_{p}(\bR^{d})$ of functions $u_{0}(x)$
can be introduced
as the space of traces on $t=0$ of 
(deterministic) functions $u $
such that
$$
u\in L_{p}(\bR_{+},H^{1}_{p}),\quad\partial u/\partial t
\in L_{p}(\bR_{+},H^{-1}_{p}),
$$
where $\bR_{+}=(0,\infty)$.
For such functions there is a 
(unique) modification denoted again
$u$ such that $u_{t}$ is a continuous $L_{p}$-valued
function on $[0,\infty)$ so that $u_{0}$ is well defined.
Any such $u_{t}$ is called an extension of $u_{0}$.

The norm in $W^{1-2/p}_{p}$ can be defined as
the infimum  of
$$
\|u\|_{ L_{p}(\bR_{+},H^{1}_{p})}
+\|\partial u/\partial t \|_{L_{p}(\bR_{+},H^{-1}_{p})}
$$
over all extensions $u_{t} $ of elements $u_{0} $.
  It is also well known
  that an equivalent norm of $u_{0}$
can be introduced as  
$$
\|u\|_{L_{p}((0,1),W^{1}_{p})},
$$
where $u=u_{t}$ is defined as the (unique)
  solution of the heat equation $\partial u_{t}(x)/
\partial t=\Delta u_{t}(x)$ with initial condition
$u_{0}(x) $.

For $s\geq0$ we introduce
$$
\tr_{\!s}\cW^{1}_{p}=L_{p}(\Omega,\cF_{s},
 W^{1-2/p}_{p}).
$$
The following auxiliary result helps understand the role
of $\tr_{\!s}\cW^{1}_{p}$. We use spaces $\cW^{1}_{p}([S,T))$
and $\bW^{1}_{p}((S,T))$,
which are introduced in the same way as $\cW^{1}_{p}(\tau)$
and $\bW^{1}_{p}(\tau)$
but the functions are only considered on $[S,T)$
and $(S,T)$, respectively.
\begin{lemma}
                                 \label{lemma 2.25.1}
Let $s\geq0$ be a fixed number and let
$u_{s}$ be an $\cF_{s}$-measurable function
with values in the set of distributions 
over $\bR^{d}$. 

(i) We have $u_{s}\in\tr_{\!s}\cW^{1}_{p}$
if and only if there
exists a $v\in\cW^{1}_{p}([s,\infty))$ satisfying
the   equation 
\begin{equation}
                                      \label{2.27.2}
\partial v/\partial t=\Delta v-v,\quad t\geq s,
\end{equation}
(which is a particular case of \eqref{11.13.1}
and is understood in the same sense)
with initial data $u_{s}$. This $v$ is unique and
satisfies
\begin{equation}
                                     \label{2.25.4}
\|v\|_{\bW^{1}_{p}((s,\infty))} \leq 
N \|u_{s}\|_{\tr_{\!s}\cW^{1}_{p}},
\quad
\|u_{s}\|_{\tr_{\!s}\cW^{1}_{p}} \leq 
N \|v\|_{\bW^{1}_{p}((s,\infty))},
\end{equation}
where the constants $N$ are independent of $s$, $u_{s}$,
and $v$.

(ii) We have $u_{s}\in\tr_{\!s}\cW^{1}_{p}$
if and only if there exists a $v\in\cW^{1}_{p}([s,s+1))$
such that $v_{s}=u_{s}$.

(iii) If such a $v$ exists and $dv_{t}=(D_{i}f^{i}_{t}
+f^{0}_{t})\,dt +g^{k}_{t}\,dw^{k}_{t}$, $t\geq s$, then
\begin{equation}
                                     \label{2.25.5}
\|u_{s}\|_{\tr_{\!s}\cW^{1}_{p}}\leq
N\big(\|v\|_{\bW^{1}_{p}((s,s+1))}+
\sum_{j=0}^{d}\|f^{j}\|_{\bL_{p}((s,s+1))}
+\|g\|_{\bL_{p}((s,s+1))}\big),
\end{equation}
where the constant $N$ is independent of $s$, $u_{s}$
and $v$.

(iv) If $s>0$ and we have a
$u\in\cW^{1}_{p}(s)$, then $u_{s}\in\tr_{\!s}\cW^{1}_{p}$
and 
$$
\|u_{s}\|_{\tr_{\!s}\cW^{1}_{p}}\leq
N\big(\|u\|_{\bW^{1}_{p}(s)}+
\sum_{j=0}^{d}\|f^{j}\|_{\bL_{p}(s)}
+\|g\|_{\bL_{p}(s)}\big),
$$
where $N$ is independent of $u$, and $f^{j}$ and $g^{k}$
are taken from \eqref{12.3.1}.

\end{lemma}

We prove this lemma in Section \ref{section 5.3.1}.

Here are our main results concerning \eqref{11.13.1}.
The following theorem 
is very close to Theorem 2.12 of \cite{Ki1}.
Important differences are the presence of the parameter
$\lambda$ in \eqref{1.4.2} and weaker assumptions on the 
coefficients of the deterministic part of the equation.

\begin{theorem}
                                    \label{theorem 12.7.1} 
  Let the above assumptions be satisfied
with    $\beta
=\beta(d,p,\delta)=\beta_{0}/3$,
where $\beta_{0} $
is the constant  from Lemma \ref{lemma 4.23.1}, and
 $\beta_{1}=\beta_{1}(d,p,\delta,
\varepsilon) >0$ taken from
 the proof of Lemma \ref{lemma 4.29.1}.
Let   $\lambda\geq0$,   let
$f^{j},g\in\bL_{p}(\tau)$, and let
$u_{0}\in\tr_{\!0}\cW^{1}_{p} $. 

(i) Then
equation \eqref{11.13.1}
for $t\leq\tau \wedge T$
has a unique solution  
$u\in\cW^{1}_{p}(\tau\wedge T)$ with initial data $u_{0}$
 and any $T\in(0,\infty)$.
Moreover, if 
$$
\lambda\geq\lambda_{0}(d,p,\delta,K,
\varepsilon,\varepsilon_{1},\varepsilon_{2})\geq1,
$$ 
 then 
equation \eqref{11.13.1}
for $t\leq\tau $
has a unique solution  
$u\in\cW^{1}_{p}(\tau )$ with initial data $u_{0}$.
 
(ii) Furthermore, if  a $v\in\cW^{1}_{p}(\infty)$ is 
defined by equation \eqref{2.27.2} with initial condition
$u_{0}$,
then the above  solution $u$   satisfies
$$
 \lambda^{1/2}\|u\|_{\bL_{p}(\tau)} 
+\|Du \|_{\bL_{p}(\tau)} 
$$
$$
  \leq N\big(\sum_{i=1}^{d}
\|f^{i}\|_{\bL_{p}(\tau)}+\|g\|_{\bL_{p}(\tau)}
+ 
\|Dv\|_{\bL_{p}(\tau)} 
\big)
$$
\begin{equation}
                                             \label{1.4.2}
+N\lambda^{-1/2} \|f^{0}\|_{\bL_{p}(\tau)}
+N\lambda^{1/2}\|v\|_{\bL_{p}(\tau)}  ,
\end{equation}
provided that $\lambda\geq\lambda_{0}$, where
the constants $N,\lambda_{0}\geq1$
depend only on $d$, $p$, $\delta$, $K$, $\varepsilon$, $\varepsilon_{1}$,
  and $\varepsilon_{2}$.

(iii) Finally, there exists a set $\Omega'\subset
\Omega$ of full probability
such that $u_{t\wedge\tau}I_{\Omega'}$ is a continuous 
$\cF_{t}$-adapted 
$L_{p}$-valued functions
of $t\in[0,\infty)$.
 
\end{theorem}

Observe that 
estimate \eqref{1.4.2} shows one of  good reasons for
writing the free term in \eqref{11.13.1}
in the form $D_{i}f^{i}+f^{0}$, because
$f^{i}$, $i=1,...,d$, and $f^{0}$
enter \eqref{1.4.2} differently.

\begin{remark}
                                           \label{remark 6.4.1}
As it follows from
our proofs, if $p=2$,
Assumptions \ref{assumption 1.2.6}
and \ref{assumption 8.9.1} are not needed
for Theorem \ref{theorem 12.7.1} to be true
and mentioning  
$\varepsilon$, $\varepsilon_{1}$, and $\varepsilon_{2}$
can be dropped in the statement.
Thus we provide a new way to
prove the classical result
on Hilbert space solvability of SPDEs
(cf., for instance, \cite{Ro}).

\end{remark}
 
 We prove Theorem \ref{theorem 12.7.1}
  in Section \ref{setion 3.2.1}
after we prepare necessary tools in 
Sections  \ref{section 5.5.3}-\ref{section 5.3.1}.
In Section \ref{section 5.5.3} we prove uniqueness
part of Theorem \ref{theorem 12.7.1} on the basis of It\^o's
formula from \cite{Trento}. Here Assumptions \ref{assumption 1.2.6}
and \ref{assumption 8.9.1} are not used.
In Section \ref{section 5.3.4} we treat the case of the heat equation
with random right-hand side and present a simplified version
of the corresponding result from 
\cite{Kr99}. In Section \ref{section 5.3.1} we prove
an auxiliary existence theorem and derive some
a priori estimates.

 Here is a result about continuous dependence
of solutions on the data.

\begin{theorem}
                                    \label{theorem 1.4.1}
 
Assume that for each $n=1,2,...$
we are given   functions $a^{ij}_{nt}$, $a^{i}_{nt}$, $
b^{i}_{nt}$, $c_{nt}$, $\sigma^{ik}_{nt}$, $\nu^{k}_{nt}$, 
$f^{j}_{nt}$, $g^{k}_{nt}$, and $u_{n0} $
having the same meaning as the original ones
 and satisfying the same assumptions
as those imposed on the original ones in
 Theorem \ref{theorem 12.7.1} 
(with the same $\delta,K,\beta,...$).
Assume that for 
$i,j=1,...,d$ and
almost all $(\omega,t,x)$ we have
$$
(a^{ij}_{nt},a^{i}_{nt},b^{i}_{nt},c_{nt})\to
(a^{ ij}_{t},a^{ i}_{t},b^{ i}_{t},c _{t}),
$$
$$
|\sigma^{ i\cdot}_{nt}-\sigma^{ i\cdot}_{t}|_{\ell_{2}}+
|\nu _{nt}-\nu _{t}|_{\ell_{2}}\to0, 
$$
as $n\to\infty$. Also  assume that 
$$
\sum_{j=0}^{d}(\|f^{ j}_{n}-f^{j}\|_{\bL_{p}(\tau)}
+
\|g_{n }-g\|_{\bL_{p}(\tau)}+
\|u_{n0} -u_{0}\|_{\tr_{\!0}\cW^{1}_{p} } \to0
$$
as $n\to\infty$. Take $\lambda\geq\lambda_{0}$, take
 the function $u$ from 
Theorem \ref{theorem 12.7.1} and let $u_{n}
\in\cW^{1}_{p}(\tau)$ be the unique solutions
of equations \eqref{11.13.1} for $t\leq\tau$ constructed from
$a^{ij}_{nt}$, $a^{i}_{nt}$, $
b^{i}_{nt}$, $c_{nt}$, $\sigma^{ik}_{nt}$, $\nu^{k}_{nt}$, 
$f^{j}_{nt}$, and $g^{k}_{nt}$ and having initial
values  $u_{n0}  $. 

Then, as $n
\to\infty$, we have
$\|u_{n}-u\|_{\bW^{1}_{p}(\tau )}\to0$ and for any
finite $T\in[0,\infty)$
\begin{equation}
                          \label{2.21.1}
E\sup_{t\leq\tau\wedge T}
\|u_{n t}-u_{t}\|_{L_{p}}^{p}\to0.
\end{equation}

\end{theorem}

Proof. Set $v_{nt}=u_{nt}-u_{t}$. Then
$$
dv_{nt}=(L_{nt}v_{nt}-\lambda v_{nt}
+D_{i}\tilde{f}^{ i}_{nt}+\tilde{f}^{ 0}_{nt})\,dt
+(\Lambda^{ k}_{nt}v_{nt}
+\tilde{g}^{ k}_{nt})\,dw^{k}_{t},
$$
where $L_{nt}$ and $\Lambda^{ k}_{nt}$ are the operators
constructed from $a^{ij}_{nt}$, $a^{i}_{nt}$, $b^{i}_{nt}$,
$c_{nt}$ and $\sigma^{ik}_{nt}$, 
$\nu^{k}_{nt}$, respectively,
and
$$
\tilde{f}^{ i}_{nt}=f^{i}_{nt}-f^{i}_{t}
+(a^{ji}_{nt}-a^{ji}_{t})D_{j}u_{t}+
(a^{ i}_{nt}-a^{i}_{t})u_{t},
$$
$$
\tilde{f}^{ 0}_{nt}=f^{ 0}_{nt}-f^{0}_{t}
+(b^{i}_{nt}-b^{i}_{t})D_{i}u_{t}
+(c_{nt}-c_{t})u_{t},
$$
$$
\tilde{g}^{ k}_{nt}=g^{k}_{nt}-g^{ k}_{t}
+(\sigma^{ik}_{nt}-\sigma^{ik}_{t})D_{i}u_{t}
+(\nu^{ k}_{nt}-\nu^{k}_{t})u_{t}.
$$

By Theorem \ref{theorem 12.7.1}
we know that $u\in\bW^{1}_{p}(\tau)$. This along with 
our assumptions and the dominated
convergence theorem implies that
$$
\sum_{j=0}^{d}\|\tilde{f}^{ j}_{n}\|_{\bL_{p}(\tau)}+
\|\tilde{g}_{n}\|_{\bL_{p}(\tau)}\to 0
$$
as $n\to\infty$. After that
by applying \eqref{1.4.2} to $v_{nt}$ 
we immediately see that
$\|u_{n}-u\|_{\bW^{1}_{p}(\tau )}\to0$.

Assertion \eqref{2.21.1} is, actually, a simple corollary
of the above. Indeed,
 by introducing $\hat{f}_{n }^{j}$
and $\hat{g}^{ k}_{n}$ in an obvious way, we can write
\begin{equation}
                                   \label{8,2,24,1}
dv_{nt}=(D_{i}\hat{f}^{i}_{nt}+\hat{f}^{0}_{nt}
)\,dt+\hat{g}^{ k}_{nt}\,dw^{k}_{t},
\end{equation}
and 
$$
\sum_{j=1}^{d}\|\hat{f}^{j}_{n}\|_{\bL_{p}(\tau )}
+\|\hat{g}_{n}\|_{\bL_{p}(\tau)}\to0.
$$
It is standard (see, for instance, our
Theorem \ref{theorem 12.3.1}) to derive from here the estimate
$$
E\sup_{t\leq\tau\wedge T}
\|u_{nt}-u_{t}\|_{L_{p}}^{p}\leq N\big(
\sum_{j=1}^{d}\|\hat{f}^{j}_{n}\|_{\bL_{p}(\tau\wedge T)}
+\|\hat{g}_{n}\|_{\bL_{p}(\tau\wedge T)}
+E 
\|u _{n0}-u_{0}\|_{L_{p}}^{p}\big),
$$
where $N$ is independent of $n$. It is also well known
that $W^{1-2/p}_{p}\subset L_{p}$, that is
$$
\|u _{n0}-u_{0}\|_{L_{p}}\leq N
\|u _{n0}-u_{0}\|_{W^{1-2/p}_{p}}.
$$
By combining all this together we obtain
\eqref{2.21.1} and the theorem is proved.

The following result could be proved 
on the basis of Theorem \ref{theorem 1.4.1}
in the same way as Corollary 5.11
of \cite{Kr99}, where the solutions are approximated
by solutions of equations with smooth coefficients and
then a stopping time techniques was used.
We give here a shorter proof
based on a different idea.

\begin{theorem}
                                    \label{theorem 1.4.2}

Let $p_{1},p_{2}\in[2,\infty)$, $p_{1}<p_{2}$, and let 
the above assumptions  be satisfied 
with $\beta\leq\beta(d,p,\delta)$ for all $p\in[p_{1},p_{2}]$
and $\beta_{1}\leq
\beta_{1}(d,p,\delta,\varepsilon)$ for all $p\in[p_{1},p_{2}]$.  
Let   $\lambda\geq0$,   and suppose that for $p\in[p_{1},p_{2}]$
we have
$f^{j},g\in\bL_{p}(\tau)$, and 
$u_{0}\in\tr_{\!0}\cW^{1}_{p} $.

 Then the solutions corresponding
to $p=p_{1}$ and $p=p_{2}$ coincide, that is, there
is a unique solution 
$u\in \cW^{1}_{p_{1}}(\tau)
\cap\cW^{1}_{p_{2}}(\tau)$
of equation \eqref{11.13.1} with initial data $u _{0}$.

\end{theorem}

Proof. Obviously, it suffices to concentrate on bounded
$\tau$. As   is explained above
in that case we may assume that $\lambda$ is as large 
as we like. We take it so large that one could use
  assertion (ii) of Theorem \ref{theorem 12.7.1}
 with any $p\in[p_{1},p_{2}]$. 

Denote by $u$ the solution corresponding to $p=p_{2}$
and observe that, owing to uniqueness of solutions
in $\cW^{1}_{p_{1}}(\tau)$, we need only show that
$u\in\cW^{1}_{p_{1}}(\tau)$.

Take a $\zeta\in C^{\infty}_{0}$ such that $\zeta(0)=1$,
set $\zeta_{n}(x)=\zeta(x/n)$,
and notice that $u^{ n}:=u \zeta_{n}$ satisfies
$$
du^{ n}_{t}=(L_{t}u^{ n}_{t})
-\lambda u^{ n}_{t}+D_{i}f^{i}_{nt}+f^{ 0}_{nt})\,dt
+(\Lambda^{k}_{t}u^{ n}_{t}+g^{k}_{nt})\,dw^{k}_{t},
$$
where
$$
f^{i}_{nt}=f^{i}_{t}\zeta_{n}-ua^{ji}_{t}
D_{j}\zeta_{n},\quad 
i\geq1,
$$
$$
f^{ 0}_{nt}=f^{0}_{t}\zeta_{n}-f^{i}_{t}D_{i}\zeta_{n}
-(a^{ij}_{t}D_{i}u_{t}+a^{j}_{t}u)D_{j}\zeta_{n}
-b^{i}_{t}u_{t}D_{i}\zeta_{n},
$$
$$
g^{k}_{nt}=g^{k}_{t}\zeta_{n}-\sigma^{ik}_{t}
u_{t}D_{i}\zeta_{n}.
$$
 It follows that for $p_{1}\leq p\leq p_{2}$ we have
\begin{equation}
                                       \label{2.24.3}
\|u^{ n}\|_{\bW^{1}_{p}(\tau)}
\leq N\big(\sum_{i=0}^{d}\|f^{ i}_{n}\|_{\bL_{p}(\tau)}
+\|g_{n}\|_{\bL_{p}(\tau)}+
\|u _{0}\zeta_{n}\|_{\tr_{\!0}\cW^{1}_{p}}\big).
\end{equation}
One knows that with constants $N$ independent of $n$
$$
\|u _{0}\zeta_{n}\|_{\tr_{\!0}\cW^{1}_{p}}
\leq N\big(\|u _{0}\zeta_{n}\|_{\tr_{\!0}\cW^{1}_{p_{1}}}+
\|u _{0}\zeta_{n}\|_{\tr_{\!0}\cW^{1}_{p_{2}}})\leq
N\big(\|u _{0} \|_{\tr_{\!0}\cW^{1}_{p_{1}}}
+\|u _{0} \|_{\tr_{\!0}\cW^{1}_{p_{2}}}).
$$
Similarly, and by H\"older's inequality
$$
\|f^{ i}_{n}\|_{\bL_{p}(\tau)}
\leq N+N\|uD\zeta_{n}\|_{\bL_{p}(\tau)}
\leq N+ \|u\|_{\bL_{p_{2}}(\tau)}\|D\zeta_{n}
\|_{\bL_{q}(\tau)},
$$
where
$$
q=\frac{pp_{2}}{p_{2}-p}.
$$
Similar estimates are available for other terms in
the right-hand side of \eqref{2.24.3}. Since
$$
\|D\zeta_{n}
\|_{\bL_{q}(\tau)}=Nn^{-1+(p_{2}-p)d/(p_{2}p)}\to0
$$
as $n\to\infty$ if
\begin{equation}
                                       \label{2.24.4}
\frac{1}{p}-\frac{1}{p_{2}}<\frac{1}{d},
\end{equation}
estimate \eqref{2.24.3} implies that $u 
\in\cW^{1}_{p}(\tau)$.

Thus knowing that $u 
\in\cW^{1}_{p_{2}}(\tau)$ allowed us to conclude
that $u 
\in\cW^{1}_{p}(\tau)$ as long as $p\in[p_{1},p_{2}]$
and \eqref{2.24.4} holds. We can now replace $p_{2}$ 
with a smaller $p$ and keep going in the same way
each time increasing $1/p$ by the same amount until $p$
  reaches $p_{1}$. Then we get that 
$u 
\in\cW^{1}_{p_{1}}(\tau)$. The theorem is proved.

In many situation the following maximum principle
is useful.

\begin{theorem}
                                    \label{theorem 1.4.3}
Let 
the above assumptions  be satisfied 
with $\beta\leq\beta(d,q,\delta)$ for all $q\in[2,p ]$
and $\beta_{1}\leq
\beta_{1}(d,q,\delta,\varepsilon)$ for all $q\in[2,p ]$. 
Let   $\lambda\geq0$ and
$f^{0}\in\bL_{p}(\tau)$, 
$u_{0}\in\tr_{\!0}\cW^{1}_{p} $,
$f^{i}=0$, $i=1,...,d$, $g=0$
be such that $u _{0}\geq0$
and $f^{0}\geq0$.
Then for the solution $u$ almost surely
we have $u_{t}\geq0$ for all finite $t\leq\tau$.
 
\end{theorem}

Proof. If $p=2$ the result is proved in
\cite{Kr07.1}. For general $p\geq2$ take the same
function $\zeta_{n}$ as in the preceding proof,
introduce $f^{n i}=f^{i}\zeta_{n}$,
$g^{ k}_{n}=0$, and call
$u^{n}$ the solution of 
\eqref{11.13.1} with so modified free terms
and the  
initial data $u _{0}\zeta_{n}$. By Theorem \ref{theorem 1.4.2}
we have $u^{n}\in\cW^{1}_{p}(\tau)\cap\cW^{1}_{2}(\tau)$.
By the above, $u^{n}\geq0$ and it only remains
to use Theorem \ref{theorem 1.4.1}. The theorem
is proved.

\mysection{It\^o's formula and uniqueness}
                                                  \label{section 5.5.3}

The following two ``standard'' results are taken from \cite{Trento}.
 
\begin{theorem}
                                     \label{theorem 12.3.1}
  Let
$u\in\cW^{1}_{p}(\tau)$, 
$f^{j}\in\bL_{p}(\tau)$,
  $g=(g^{k})\in\bL_{p}(\tau)$ and assume that
\eqref{12.3.1} holds
for $t\leq\tau$ in the sense of generalized functions.
Then there is a set $\Omega'\subset\Omega$ of full probability
such that

(i) $u_{t\wedge\tau}I_{\Omega'}$
is a continuous $L_{p}$-valued $\cF_{t}$-adapted function on
$[0,\infty)$;

(ii) for all $t\in[0,\infty)$ and $\omega\in\Omega'$ It\^o's formula holds:
$$
\int_{\bR^{d}}|u_{t\wedge\tau}|^{p}\,dx
=\int_{\bR^{d}}|u_{0}|^{p}\,dx
+p
\int_{0}^{t\wedge\tau }\int_{\bR^{d}}|u _{s}|^{p-2}
u _{s}
g^{k }_{s}\,dx\,dw^{k}_{s}
$$
$$+
\int_{0}^{t\wedge\tau }
\big( \int_{\bR^{d}}\big[p|u_{t}|^{p-2}u_{t}f^{0}_{t}
-p(p-1)|u_{t}|^{p-2}f^{i}_{t}D_{i}u_{t}
$$
\begin{equation}
                                            \label{4.19.5}
+(1/2)p(p-1)|u_{t}|^{p-2}|g_{t}|_{\ell_{2}}^{2}
\big]\,dx\big)\,dt.
\end{equation}

Furthermore,  for any $T\in[0,\infty)$ 
$$
 E\sup_{t\leq\tau\wedge T}
\|u_{t}\|^{p}_{L_{p}}\leq  2E\|u_{0}\|^{p}_{L_{p}}+
NT^{p-1}\|f^{0}\|^{p}_{\bL_{p}(\tau)}
$$
\begin{equation}
                                         \label{4.11.5}
+NT^{(p-2)/2}(\sum_{i=1}^{d}\|f^{i}\|^{p}_{\bL_{p}(\tau)}
+\|g\|^{p}_{\bL_{p}(\tau)}+\|Du\|^{p}_{\bL_{p}(\tau)}) ,
\end{equation}
where $N=N(d,p)$.  
\end{theorem}

Here is an ``energy" estimate.

\begin{corollary}
                               \label{corollary 4.19.1}
Under the conditions of Theorem \ref{theorem 12.3.1}
assume that $\tau<\infty$ (a.s.). Then
$$
E\int_{\bR^{d}}|u_{0}| ^{p}\,dx+E\int_{0}^{\tau}
\big( \int_{\bR^{d}}\big[p|u_{t}|^{p-2}u_{t}f^{0}_{t}
-p(p-1)|u_{t}|^{p-2}f^{i}_{t}D_{i}u_{t}
$$
\begin{equation}
                                       \label{12.3.2}
+(1/2)p(p-1)|u_{t}|^{p-2}|g_{t}|_{\ell_{2}}^{2}
\big]\,dx\big)\,dt 
\geq EI_{\tau<\infty}\int_{\bR^{d}}|u_{\tau}| ^{p}\,dx.
\end{equation}
Furthermore, if $\tau$ is bounded then there is an equality
instead of inequality in \eqref{12.3.2}.
\end{corollary}

The next result implies, in particular, uniqueness
in Theorem \ref{theorem 12.7.1}.

\begin{lemma}
                                \label{lemma 12.4.1}
Under Assumption \ref{assumption 1.2.1}
there exist $\lambda_{0}\geq0$ and $N$ depending only on
$d,p,K $, and $\delta$
such that, for any 
strictly positive $\lambda\geq\lambda_{0}$
and any solution $u\in\cW^{1}_{p,0}(\tau)$
 of \eqref{11.13.1} for $t\leq\tau$,
we have
\begin{equation}
                                       \label{12.3.3}
\lambda\|u\|_{\bL_{p}(\tau)}\leq N
\lambda^{1/2}\big(\sum_{j=1}^{d}
\|f^{j}\|_{\bL_{p}(\tau)}+\|g\|_{\bL_{p}(\tau)}
\big)+N\|f^{0}\|_{\bL_{p}(\tau)}.
\end{equation}
Furthermore, if $a^{i}=b^{i} =\nu^{k}\equiv0$,
then one can take $\lambda_{0}=0$.
\end{lemma}

Proof. We may assume that $f^{j}\in\bL_{p}(\tau)$,
  $g=(g^{k})\in\bL_{p}(\tau)$, since otherwise
the right-hand side of \eqref{12.3.3} is infinite.

If \eqref{12.3.3} is true for $\tau\wedge T$
in place of $\tau$
and any $T\in(0,\infty)$, then it is obviously
also true as is. Therefore, we may assume that
$\tau$ is finite. An advantage of this assumption
is that   we can use Corollary \ref{corollary 4.19.1}.
Write \eqref{12.3.2}
with $\hat{f}^{i}_{t}$, $\hat{f}^{0}_{t} $, and
$\hat{g}^{k}_{t} $ in place of $f^{i}_{t}$,
$f^{0}_{t}$, and $g^{k}_{t}$, respectively,
where
$$
\hat{f}^{i}_{t}=
a^{ji}_{t}D_{j}u_{t}+ a^{i}_{t}u_{t}+f^{i}_{t},
$$
$$
\hat{f}^{0}_{t}=
b^{i}_{t}D_{i}u_{t}+(c_{t}-\lambda)u_{t}+f^{0}_{t},\quad
\hat{g}^{k}_{t}=\sigma^{ik}_{t}D_{i}u_{t}
+\nu^{k}_{t}u_{t}+g^{k}_{t}.
$$
Then observe that inequalities like
$(a+b)^{2}\leq(1+\varepsilon)a^{2}+(1+\varepsilon^{-1})
b^{2}$ show that for any $\varepsilon\in(0,1]$ we have
$$
 |\hat{g}_{t}|_{\ell_{2}}^{2}\leq (1+\varepsilon)
\big|\sum_{i=1}^{d}\sigma^{i\cdot}_{t}D_{i}u_{t}\big|^{2}
_{\ell_{2}}
+2\varepsilon^{-1}|\nu_{t}u_{t}+g_{t}|^{2}_{\ell_{2}}
$$
$$
\leq2(1+\varepsilon)\alpha^{ij}_{t}(D_{i}u_{t})D_{j}u_{t}
+N\varepsilon^{-1}(|u_{t}|^{2}+|g_{t}|_{\ell_{2}}^{2}).
$$
Owing to \eqref{1.3.2}, for $\varepsilon=\varepsilon
(\delta)>0$ small enough
$$
I_{t}:= (1/2) |u_{t}|^{p-2}|\hat{g}_{t}|_{\ell_{2}}^{2} 
-|u_{t}|^{p-2}\hat{f}^{i}_{t}D_{i}u_{t}+(p-1)^{-1}
|u_{t}|^{p-2}u_{t}b^{i}_{t}D_{i}u_{t}
$$
\begin{equation}
                                               \label{1.3.4}
\leq-(\delta/2)|u_{t}|^{p-2}|Du_{t}|^{2}+
N|u_{t}|^{p-2}(|u_{t}|^{2}+|g_{t}|_{\ell_{2}}^{2}
+|Du_{t}|\,|u_{t}|+|Du_{t}|\sum_{i=1}^{d}|f^{i}_{t}|).
\end{equation}
Next we use that for any $\gamma>0$
$$
|u_{t}|^{p-1}|Du_{t}|=(|u_{t}|^{(p-2)/2}|Du_{t}|)
|u_{t}|^{p/2}\leq\gamma|u_{t}|^{p-2}|Du_{t}|^{2}
+\gamma^{-1} |u_{t}|^{p},
$$
$$
|u_{t}|^{p-2}|Du_{t}|\,|f^{i}_{t}|
\leq\gamma|u_{t}|^{p-2}|Du_{t}|^{2}+\gamma^{-1}
|u_{t}|^{p-2}|f^{i}_{t}|^{2},
$$
and by choosing $\gamma$ appropriately find from
\eqref{1.3.4} that
\begin{equation}
                                               \label{1.3.5}
I_{t}\leq N|u_{t}|^{p}+N|u_{t}|^{p-2}
\big(\sum_{i=1}^{d}|f^{i}_{t}|^{2}+|g_{t}|_{\ell_{2}}^{2}
\big).
\end{equation}

After that   H\"older's inequality and
\eqref{12.3.2}, where the right-hand side
is nonnegative, immediately lead to 
$$
( \lambda-N_{1})\|u\|_{\bL_{p}(\tau)}^{p}
\leq N\|u\|_{\bL_{p}(\tau)}^{p-2}\big(\sum_{i=1}^{d}
\|f^{i}\|_{\bL_{p}(\tau)}^{2}+\|g\|_{\bL_{p}(\tau)}^{2}\big)
+N\|u\|_{L_{p}(\tau)}^{p-1}\|f^{0}\|_{L_{p}(\tau)}.
$$
Furthermore, simple inspection of the above argument shows
that, if $a^{i}=b^{i} =\nu^{k}\equiv0$,
then the terms with $|u_{t}|^{2}$ and $|u_{t}|\,|Du_{t}|$
in \eqref{1.3.4} and the term with $|u_{t}|^{p}$ in 
\eqref{1.3.5} disappear, so that 
we can take $N_{1}=0$ in this case (recall that $c\leq0$).
  Generally,
for $\lambda\geq2N_{1}$ we have $\lambda-N_{1}\geq(1/2)\lambda$
and
$$
\bar{U}^{p}\leq N\bar{U}^{p-2}\bar{G}^{2}+N\bar{U}^{p-1}
 \bar{F},
$$
where
$$
\bar{U}=\lambda\|u\|_{\bL_{p}(\tau)},\quad
\bar{G}=\lambda^{1/2}(\|f^{i}\|_{\bL_{p}(\tau)} +
\|g\|_{\bL_{p}(\tau)}),\quad\bar{F}=\|f^{0}\|_{L_{p}(\tau)}.
$$

It follows that $\bar{U}\leq N(\bar{G}+\bar{F})$,
which is \eqref{12.3.3} and the lemma is proved.

\mysection{Case of the heat equation}
                                                  \label{section 5.3.4}

To move further we need the following analytic
fact established in \cite{Kr94} (see also
\cite{Kr06} for a complete proof).
\begin{lemma}
                                              \label{lemma 4.11.1}
Denote by $T_{t}$ the heat semigroup in $\bR^{d}$
and let $p\geq2$,
 $-\infty\leq a<b\leq\infty$, $g\in L_{p}((a,b)\times\bR^{d},
\ell_{2})$. Then
$$
\int_{\bR^{d}}\int_{a}^{b}
\big[\int_{a}^{t}
|DT_{t-s}g_{s}(x)|_{\ell_{2}}^{2}\,ds\big]\,dtdx
\leq N(d,p)\int_{\bR^{d}}\int_{a}^{b}
|g_{t}(x)|_{\ell_{2}}^{2}\,dtdx.
$$

\end{lemma}

In this section we deal with the following
   model equation
\begin{equation}
                                           \label{4.9.1}
du_{t}=\Delta u_{t}\,dt+g^{k}_{t}\,dw^{k}_{t}.
\end{equation}
\begin{lemma}
                                                \label{lemma 4.9.1}
Assume that $\tau\leq T$, where the
  constant $T\in[0,\infty)$.
Then for any $g=(g^{1},g^{2},...)\in\bL_{p}(\tau)$
there exists a unique $u\in\cW^{1}_{p,0}(\tau)$
satisfying \eqref{4.9.1} for $t\leq\tau$.
Furthermore, for this solution we have
\begin{equation}
                                       \label{4.10.1}
E\sup_{t\leq\tau}\|u_{t}\|
^{p}_{ L_{p} }\leq N(d,p)T^{(p-2)/2}\|g\|^{p}_{\bL_{p}(\tau)},
\end{equation}
\begin{equation}
                                       \label{4.10.2}
\|Du\|_{\bL_{p}(\tau)}\leq N(d,p) \|g\|_{\bL_{p}(\tau)}.
\end{equation}
 \end{lemma}

Proof. By replacing the unknown function
$u_{t}$ with $v_{t}e^{\lambda t}$ we see that
$v_{t}$ satisfies
$$
dv_{t}=(\Delta v_{t}-\lambda v)\,dt+e^{-\lambda t}
g^{k}_{t}\,dw^{k}_{t}.
$$
Since $\tau$ is bounded the inclusions
$u\in\cW^{1}_{p,0}(\tau)$ and $v\in\cW^{1}_{p,0}(\tau)$
are equivalent and our assertion about uniqueness
follows from   Lemma \ref{lemma 12.4.1}. 

In the proof of existence we borrow part of the proof
of Lemma 4.1 of \cite{Kr99}. As we have pointed out in the
Introduction, the beginning of the theory of divergence and 
nondivergence type equations is the same. The only difference
 with that proof is that here we take $f\equiv 0$.

We take an integer $m\geq1$, some bounded stopping times
$\tau_{0}\leq\tau_{1}\leq...\leq\tau_{m}\leq T$ and some 
(nonrandom) functions
$g^{ij} \in C^{\infty}_{0}$, $i,j=1,...,m$. Then we define
$$
g^{k}_{t}(x)=\sum_{i=1}^{m}g^{ik}(x)I_{(\tau_{i-1},\tau_{i}]}(t),
$$
$$
v_{t}( x)=\sum_{k=1}^{m}
\int_{0}^{t}g^{k}_{s}(x)\,dw^{k}_{s}
=\sum_{i,k=1}^{m}g^{ik}(x)(w^{k}_{t\wedge\tau_{i}}
-w^{k}_{t\wedge\tau_{i-1}}),\quad t\geq 0.
$$
Obviously, for any $\omega$, the function $v_{t}(x)$
is continuous and bounded in $(t,x)$ along with any derivative
in $x$. Furthermore, the function and its derivative in $x$
are H\"older $1/3$ continuous in $t$ uniformly with respect to $x$
(for almost any $\omega$). Also $v_{t}(x)$ has compact support in $x$.

These properties of $v_{t}(x)$ imply that for any $\omega$
there exists a unique classical solution of the heat equation
$$
\frac{\partial}{\partial t}\bar{u}_{t}=\Delta\bar{u}_{t}+\Delta
v_{t},\quad t>0,
$$
with zero initial data. Furthermore, 
\begin{equation}
                                                  \label{4.11.2}
\bar{u}_{t}(x)=\int_{0}^{t}T_{t-s}\Delta v_{s}(x)\,ds.
\end{equation}
This formula shows, in particular, that $\bar{u}_{t}(x)$
is $\cF_{t}$-adapted. 
Adding the fact that $\bar{u}_{t}$ is continuous
in $t$ proves that $\bar{u}_{t}(x)$ is predictable.
The same holds for
$$
(\bar{u}_{t},\phi)=\int_{0}^{t}(T_{t-s}\Delta v_{s},\phi)\,ds
$$
with any $\phi\in C^{\infty}_{0}$. 
The following corollary of Minkowski's inequality
\begin{equation}
                                                  \label{4.11.3}
\|\bar{u}_{t}\|_{L_{p}}\leq
\int_{0}^{t}\|\Delta v_{s}\|_{L_{p}}\,ds
\end{equation}
shows that $\bar{u}_{t}$ is $L_{p}$-valued. Since
$(\bar{u}_{t},\phi)$ is predictable for any
$\phi\in C^{\infty}_{0}$, $\bar{u}_{t}$ is weakly
and hence strongly   predictable as an $L_{p}$-valued
process. 

One can differentiate \eqref{4.11.2} with respect to $x$
as many times as one wants and get similar statements 
about the derivatives of $\bar{u}_{t}$. In particular,
\eqref{4.11.3} implies that
for any multi-index
$\alpha$
$$
E\int_{0}^{T}\int_{\bR^{d}}
|D^{\alpha}\bar{u}_{t}|^{p}\,dxdt
\leq T^{p}E\int_{0}^{T}\int_{\bR^{d}}
|D^{\alpha}\Delta v_{t}|^{p}\,dxdt<\infty,
$$
so that $\bar{u}_{t}\in\cW^{1}_{p,0}(T)$.

Now, it is easily seen that
$$
u_{t}(x):=\bar{u}_{t}(x)+v_{t}(x)
$$
satisfies \eqref{4.9.1} pointwisely
and by the above 
$u_{t}\in\cW^{1}_{p,0}(T)$. The (deterministic)
  Fubini's theorem also shows that $u_{t}$
satisfies \eqref{4.9.1} in the sense of distributions.

Next, we use the same simple transformation as in the proof
of Lemma 4.1 of \cite{Kr99} and conclude that for any $t$ and $x$
almost surely
$$
Du_{t}(x)=\sum_{k=1}^{m}\int_{0}^{t}
T_{t-s}Dg^{k}_{s}(x)\,dw^{k}_{s}.
$$
Hence by Burkholder-Davis-Gundy inequality
$$
E|Du_{t}(x)|^{p}\leq NE\big[
\int_{0}^{t}|
T_{t-s}Dg _{s}(x)|^{2}_{\ell_{2}}\,ds\big]^{p/2},
$$
which along with Lemma \ref{lemma 4.11.1} proves
\eqref{4.10.2} for our particular $g$.
Theorem \ref{theorem 12.3.1} shows that  
\eqref{4.10.1} follows from \eqref{4.10.2}
and \eqref{4.9.1}.

The rest is trivial since the set of $g$'s
like the one above
is dense in $\bL_{p}(T)$ by Theorem 3.10
of \cite{Kr99}. The lemma is proved.

Next we introduce the parameter $\lambda$ into
\eqref{4.9.1}.
 
\begin{lemma}
                                                \label{lemma 4.9.01}
Assume that $\tau\leq T$, where the
  constant $T\in[0,\infty)$.
Let $\lambda>0$.
Then for any $g=(g^{1},g^{2},...)\in\bL_{p}(\tau)$
there exists a unique $u\in\cW^{1}_{p,0}(\tau)$
satisfying 
\begin{equation}
                                           \label{4.21.3} 
du_{t}=(\Delta u_{t}-\lambda u_{t})\,dt+g^{k}_{t}\,dw^{k}_{t}.
\end{equation} for $t\leq\tau$.
Furthermore, for this solution we have
\begin{equation}
                                       \label{4.21.7}
\lambda^{p/2}\|u\|^{p}_{\bL_{p}(\tau)}
\leq N(d,p) \|g\|^{p}_{\bL_{p}(\tau)},
\end{equation}
\begin{equation}
                                       \label{4.21.8} 
\|Du\|_{\bL_{p}(\tau)}\leq N(d,p) \|g\|_{\bL_{p}(\tau)}.
\end{equation}
 \end{lemma}

Proof. Uniqueness and estimate \eqref{4.21.7}
follow from Lemma \ref{lemma 12.4.1}. The existence
immediately follows from Lemma \ref{lemma 4.9.1}
and the result of transformation described in the beginning
of its proof. To establish \eqref{4.21.8} consider the heat
equation
\begin{equation}
                                       \label{4.21.5}
\frac{\partial}{\partial t}v_{t}=\Delta v_{t}-\lambda u_{t}.
\end{equation}
Since $u\in\bL_{p}(\tau)$, 
for almost any $\omega$ we have $u\in L_{p}((0,\tau)
\times\bR^{d})$ and by by a classical result 
(see, for instance, \cite{Kr08_2})
 for almost any $\omega$ equation \eqref{4.21.5}
with zero initial data has a unique solution
in the class of functions such that along with derivatives in $x$
up to the second order they belong to 
$  L_{p}((0,\tau)
\times\bR^{d})$. Furthermore,
$$
\|D^{2}v\|^{p}_{L_{p}((0,\tau)
\times\bR^{d})}+
\lambda^{p/2}\|D v\|^{p}_{L_{p}((0,\tau)
\times\bR^{d})}
$$
\begin{equation}
                                       \label{4.21.6}
+\lambda^{p}\| v\|^{p}_{L_{p}((0,\tau)
\times\bR^{d})}\leq N\|\lambda u\|^{p}_{L_{p}((0,\tau)
\times\bR^{d})}.
\end{equation}
The solution $v_{t}$ can be given by an integral formula,
which implies that $v_{t}$ is 
$\cF_{t}$-adapted. It is also continuous as an $L_{p}$-valued
process, hence, is
a predictable $L_{p}$-valued process.
Taking expectations
of  both parts of \eqref{4.21.6}  
  shows that $v\in\cW^{1}_{p}(\tau)$.

Now observe that
$$
d(u_{t}-v_{t})=\Delta(u_{t}-v_{t})\,dt+g^{k}_{t}\,dw^{k}_{t},
$$
which by Lemma \ref{lemma 4.9.1} implies that
$$
\|D(u-v)\|^{p}_{\bL_{p}(\tau)}\leq N\|g\|^{p}_{\bL_{p}(\tau)}.
$$
Upon combining this with \eqref{4.21.6} we obtain
$$
\|Du\|^{p}_{\bL_{p}(\tau)}\leq N(\|g\|^{p}_{\bL_{p}(\tau)}
+\lambda^{p/2}\|u\|^{p}_{\bL_{p}(\tau)}),
$$
which along with \eqref{4.21.7} yields \eqref{4.21.8}.
The lemma is proved.

\mysection{A priori estimates in the general case}

                                        \label{section 5.3.1}

First we deal with the case when $\sigma=\nu=0$.  

\begin{lemma}
                                    \label{lemma 4.23.1} 
Suppose 
that    $\sigma^{ik}\equiv\nu^{k}\equiv0$. 
Also suppose that Assumptions \ref{assumption 1.2.1}
and \ref{assumption 1.2.6}  
are satisfied with $\beta\leq\beta_{0}$,
where the way to estimate the
constant $\beta_{0}(d,p,\delta)>0$ is described in the proof.
Let    $f^{j}\in\bL_{p}(\tau)$
 and $g \in\bL_{p}(\tau)$. 

Then there exist constants $\lambda_{0}\geq1$
and $N$, depending only on $d,p,\delta,K $, and  
$\varepsilon$, such that for any $\lambda\geq\lambda_{0}$
there exists
 a unique $u\in\cW^{1}_{p,0}(\tau)$
satisfying \eqref{11.13.1} for $t\leq\tau$. Furthermore, this solution
satisfies the estimate
\begin{equation}
                                             \label{11.18.3}
 \lambda^{1/2}\|u\|_{\bL_{p}(\tau)} 
+\|Du \|_{\bL_{p}(\tau)} 
\leq N\big(\sum_{i=1}^{d}
\|f^{i}\|_{\bL_{p}(\tau)}+\|g\|_{\bL_{p}(\tau)}
\big)+N\lambda^{-1/2}\|f^{0}\|_{\bL_{p}(\tau)}.
\end{equation}
  
\end{lemma}

Proof. Uniqueness and part of estimate \eqref{11.18.3}
follow from Lemma \ref{lemma 12.4.1}. 
In the rest of the proof we may assume that $\tau$ is bounded
and split our argument into two parts.

{\em Case $g^{k}\equiv0$\/}. 
First assume that the coefficients and $f^{j}$
are nonrandom.
We extend the coefficients of $L$ following the example
$a^{ij}_{t}(x)=\delta^{ij}$, $t<0$,
and extend $f^{j}_{t}$   beyond $(0,\tau)$ arbitrary
 only requiring $f^{j}\in L_{p}(\bR^{d+1})$.

Then  by Theorem  4.5 
and Remark 2.4 of \cite{Kr07} the equation
\begin{equation}
                                                      \label{4.24.3}
\frac{\partial}{\partial t}u_{t}=
L_{t}u_{t}-\lambda u_{t}
+D_{i}f^{i}_{t}+f^{0}_{t}
\end{equation} 
in $\bR^{d+1}$ has a unique solution with finite norms
$$
\|u\|_{L_{p}(\bR^{d+1})}\quad\text{and}\quad
\|Du\|_{L_{p}(\bR^{d+1})}
$$
provided that $\lambda\geq\lambda_{0}$.
By Theorem  4.4 of \cite{Kr07}
\begin{equation}
                                                     \label{4.24.1}
\lambda^{1/2}\|u\|_{L_{p}(\bR^{d+1})}+
\|Du\|_{L_{p}(\bR^{d+1})}\leq N(\sum_{i=1}^{d}\|f^{i}
\|_{L_{p}(\bR^{d+1})}
+\lambda^{-1/2}\|f^{0}\|_{L_{p}(\bR^{d+1})}).
\end{equation}  
By Theorem \ref{theorem 12.3.1} the function
$u_{t}$ is a continuous $L_{p}$-valued function.

The proof of Theorem  4.4 of \cite{Kr07} is achieved
on the basis of the a priori estimate \eqref{4.24.1}
and the method of
continuity by considering the family of equations
\begin{equation}
                                                \label{4.24.2}
\frac{\partial}{\partial t}u_{t}=
(\theta L_{t}+(1-\theta)\Delta)u_{t}-\lambda u_{t}
+D_{i}f^{i}_{t}+f^{0}_{t},
\end{equation}  
where the parameter $\theta$ changes in $[0,1]$.
We remind briefly the method of continuity because we 
want to show that certain properties of equation
\eqref{4.24.2} which we know for $\theta=0$ propagate
from $\theta=0$ to $
\theta=1$. 

We fix a $\theta_{0}\in[0,1]$ and to solve \eqref{4.24.2}
for given $f^{j}$ define a sequence of
$u^{n}\in L_{p}(\bR,W^{1}_{p})$  by solving
  the equation
$$
\frac{\partial}{\partial t}u^{n+1}_{t}=
(\theta_{0} L_{t}+(1-\theta_{0})\Delta)u^{n+1}_{t}-
\lambda u^{n+1}_{t}
$$
\begin{equation}
                                                \label{5.8.1}
+D_{i}f^{i}_{t}+f^{0}_{t}+(\theta-\theta_{0})(L_{t}-\Delta)u^{n},
\quad n\geq1,\quad u^{0}=0.
\end{equation} 
If we know that equation \eqref{4.24.2} is uniquely solvable
with $\theta_{0}$ in place of $\theta$ for arbitrary
$f^{j}\in L_{p}(\bR^{d+1})$, then the sequence $u^{n}$
is well defined. Furthermore, estimate \eqref{4.24.1}
easily shows that for $\theta$  sufficiently close
to $\theta_{0}$ the $L_{p}(\bR,W^{1}_{p})$
norm of $u^{n+1}-u^{n}$ goes to 
zero geometrically as $n\to\infty$. In this way passing to the limit
in \eqref{5.8.1} we obtain the solution of 
\eqref{4.24.2} for $\theta$ close to $\theta_{0}$. Then we can repeat
the procedure and starting from $\theta=0$ and moving
 step by step eventually reach $\theta=1$.

For $\theta=0$ we are dealing with solvability of the heat
equation which is proved by giving
the solution explicitly by means of the heat semigroup.
This representation formula 
has two important implications:

(i) For any constant $T\in\bR$,
changing $f^{j}_{t}$ for $t\geq T$ does not affect
$u_{t}$ for $t\leq T$;

(ii) If $f^{j}$ are $L_{p}(\bR^{d+1})$-valued
measurable functions of a parameter, say $\omega$
from a measurable space, say $(\Omega,\cF_{T})$, then   the solution
 $u\in L_{p}(\bR,W^{1}_{p})$, which now depends on $\omega$
 is also $\cF_{T}$-measurable.

Property (i) is obtained by inspecting the representation formula.
Property (ii) is true  because 
the mapping $L_{p}(\bR^{d+1})\ni f^{j}\to u
\in L_{p}(\bR,W^{1}_{p})$ is continuous and hence
Borel measurable. 

Obviously, both properties
propagate from $\theta=0$ to $\theta=1$ by
the above method of continuity. In particular, 
solutions of \eqref{4.24.3} 
on the time interval $(-\infty,T]$ depend
only on the values of $f^{j}_{t}$ for $t\in(-\infty,T]$.
 It follows that with the same
$\lambda$ and $N$, for any $T\in\bR$,
$$
\lambda^{1/2}\|u\|_{L_{p}((-\infty,T),L_{p})}+
\|Du\|_{L_{p}((-\infty,T),L_{p})}
$$
\begin{equation}
                                                       \label{4.24.4}
\leq N(\sum_{i=1}^{d}\|f^{i}\|_{L_{p}((-\infty,T),L_{p})}
+\lambda^{-1/2}\|f^{0}\|_{L_{p}((-\infty,T),L_{p})}.
\end{equation} 

From now on 
we allow the coefficients and $f^{j}$ to be random,
 continue $f^{j}$ as zero 
for $t<0$ and solve \eqref{4.24.3} for each $\omega$.
By \eqref{4.24.4} with $T=0$ we have that  $u_{t}=0$ 
for $t\leq0$ and it makes sense considering equation
\eqref{4.24.3} on $(0,T)$ for each $T\in(0,\infty)$
with zero initial condition. In such situation
 properties (i) and (ii) still hold.

In particular,
 if $f^{j} $ are measurable $L_{p}((0,T),L_{p})$-valued functions
of a parameter, say $\omega$
from a measurable space, say $(\Omega,\cF_{T})$, then   the solution
 $u\in L_{p}((0,T),W^{1}_{p})$ 
 is also $\cF_{T}$-measurable.
Then from the equation itself it follows that
$(u_{T},\phi)$ is $\cF_{T}$-measurable for any 
$\phi\in C^{\infty}_{0}$. Since $u_{T}$ takes values
in $L_{p}$, it is an $L_{p}$-valued $\cF_{T}$-measurable
function. 

 If $f^{i}_{t}$ are predictable $L_{p}$-valued
function, the above conclusions are valid for any
$T\in[0,\infty)$. In particular, $u_{t}$
is $\cF_{t}$-adapted as an $L_{p}$-valued function
and since it is continuous, 
$u_{t}$ is a predictable $L_{p}$-valued function.

  These properties
and the fact that \eqref{4.24.4} holds for any
$T\in(0,\infty)$ and $\omega$ prove the lemma in the particular case
under consideration.
 
{\em General case\/}. By Lemma \ref{lemma 4.9.01}
there is a unique solution $v\in\cW^{1}_{p,0}(\tau)$ of
\eqref{4.21.3}. 
Observe that
$$
(L_{t}-\Delta)v_{t}=D_{i}\hat{f}^{i}_{t}+\hat{f}^{0}_{t},
$$
where $\hat{f}^{j}_{t}$ are function of class $\bL_{p}(\tau)$
defined by
$$
\hat{f}^{j}_{t}
=(a^{ij}_{t}-\delta^{ij})D_{i}v_{t}+a^{j}_{t}v_{t},\quad j=1,...,d,
$$
$$
\hat{f}^{0}_{t}=b^{i}_{t}D_{i}v_{t}+c_{t}v_{t}.
$$

By the above there is a unique solution
$u\in\cW^{1}_{p,0}(\tau)$ of 
$$
\frac{\partial}{\partial t}u_{t}=
L_{t}u_{t}-\lambda u_{t}+(L_{t}-\Delta)v_{t} 
+D_{i}f^{i}_{t}+f^{0}_{t}.
$$
Obviously, $v_{t}+u_{t}$ is a solution
of class $\cW^{1}_{p,0}(\tau)$ of equation
\eqref{11.13.1}. By the particular case
$$
\lambda^{1/2}\|u\|_{\bL_{p}(\tau)} 
+\|Du \|_{\bL_{p}(\tau)} 
\leq N\big(\sum_{i=1}^{d}(
\|f^{i}\|_{\bL_{p}(\tau)} +\|\hat{f}^{i}\|_{\bL_{p}(\tau)})
\big)
$$
$$
+N\lambda^{-1/2}(\|f^{0}\|_{\bL_{p}(\tau)}
+\|\hat{f}^{0}\|_{\bL_{p}(\tau)})
$$
and to obtain \eqref{11.18.3} it only remains
to use the estimates of $v_{t}$ provided by Lemma
\ref{lemma 4.9.01}. The lemma is proved.

Now we allow $\sigma\ne0$.

\begin{lemma}
                                    \label{lemma 4.29.1}

(i) Suppose that Assumptions \ref{assumption 1.2.1}
is satisfied with $K=0$ and take   $\varepsilon\geq\varepsilon_{1}>0$,
 $\varepsilon_{2}\in(0,\varepsilon/4]$, $t_{0}\geq0$,
and $x_{0}\in\bR^{d}$. 

(ii) Let    $f^{j}\in\bL_{p}(\tau)$,
 $g \in\bL_{p}(\tau)$, and  $u\in\cW^{1}_{p,0}(\tau)$
be such that \eqref{11.13.1} holds for $t\leq\tau$.  
Assume that $u_{t}(x)=0$ if 
$$
(t,x)\not\in \Gamma:=(t_{0},t_{0}+\varepsilon_{1}^{2})
\times B_{\varepsilon_{2}}(x_{0}).
$$

(iii)
Assume that  the couple $(a,\sigma)$
is $(\varepsilon,\varepsilon_{1})$-regular
at $(t_{0},x_{0})$ with $\beta=\beta_{0}/3$
in \eqref{8.7.1} and \eqref{8.6.1},
where $\beta_{0} $
is the constant  from Lemma \ref{lemma 4.23.1}.      
Also assume that
$$
|\sigma^{i\cdot}_{t}(x)-
\sigma^{i\cdot}_{t}(x_{0})|_{\ell_{2}}\leq\beta_{1},\quad 
(a^{jk}_{t}(y)  
-  \alpha^{jk}_{t}(x_{0})) \xi^{j}
\xi^{k}\geq\delta|\xi|^{2}
$$
for all values of indices and arguments 
such that $(t,x)\in\Gamma$
and $(t,y)\in Q_{\varepsilon }(t_{0},x_{0})$, where $\beta_{1}=
\beta_{1}(d,\delta,p,\varepsilon)>0$  is a constant
an estimate from below for which can be obtained   
from the proof.

Then there exist constants $\lambda_{0}\geq1$
and $N$, depending only on $d$, $p$, $\delta$,     
  and $\varepsilon $, such that 
estimate \eqref{11.18.3} holds
provided that $\lambda\geq\lambda_{0}$.

\end{lemma}

Proof.  Without loss of generality we may and will assume
that $x_{0}=0$.
Also we modify, if necessary, $a$ and $\sigma$ in such a way that
$\sigma^{ik}_{t}(x)=0$ if $t\not\in(t_{0},t_{0}
+\varepsilon_{1}^{2})$, and $a^{ij}_{t}(x)=\delta^{-1}\delta^{ij}$
if $t\not\in(t_{0},t_{0}
+\varepsilon_{1}^{2})$. Obviously, under this modification
assumption (iii) is preserved and equation \eqref{11.13.1}
remains unaffected due to assumption (ii).
 The rest of the proof we split into two cases.

{\em Case $\sigma^{ik}_{t}(x)=
\sigma^{ik}_{t}(0)$ for $|x|\leq\varepsilon_{2}$ and $t\geq0$}.
 We want to apply Lemma \ref{lemma 4.23.1} and for that,
even if $\sigma\equiv0$,
we need $a^{ij}$ to satisfy at least the condition  
$\text{osc\,}(a^{ij},Q)\leq\beta$
for {\em all\/}
 $Q\in\bQ$ with $\rho(Q)\leq\varepsilon$. To achieve this we
modify $a^{ij}_{t}(x)$ for $|x|\geq\varepsilon/4$ using the fact that
such   modifications have no effect on \eqref{11.13.1} since
$u_{t}(x)=0$ for $|x|\geq\varepsilon_{2}$
and $\varepsilon_{2}\leq\varepsilon/4$.

Take 
a   $\xi\in C^{\infty}_{0}(\bR^{d})$
with   support lying in the ball of
radius $\varepsilon/2$ centered at the origin
and such that $\xi(x)=1$ for $|x|\leq\varepsilon/4$ and $0\leq\xi\leq1$.
Set
$$
\hat{a}^{ij}_{t}:=\xi a^{ij}_{t}+\delta^{-1}(1-\xi)\delta^{ij}.
$$
We can use $\hat{a}$ in place of $a$
in \eqref{11.13.1}. It follows
by Lemma 4.7 of \cite{Kr99} 
(It\^o-Wentzell formula) that
the function $v_{t}(x ):=
u_{t}(x+x_{t})$ satisfies   the equation
\begin{equation}
                                          \label{7.5.1}
dv_{t}(x)=(\bar{L}_{t}v_{t}(x)-\lambda v_{t}
+D_{i}\bar{f}^{i}_{t}+\bar{f}^{0}_{t})\,dt
+ \bar{g}^{k}_{t}(x+x_{t}) \,dw^{k}_{t},
\end{equation}
where     
$$
\bar{L}_{t}\phi=D_{j}(\bar{a}^{ij}_{t}D_{i}\phi),\quad
 \bar{a}^{ij}_{t}(x)= 
\hat{a}^{ij}_{t}(x+x_{t})-\alpha^{ij}_{t}(0),
$$
$$
\bar{f}^{i}_{t}(x):=f^{i}_{t}(x+x_{t})-\sigma_{t}^{ik}(0)g^{k}_{t}
(x+x_{t}),\quad i=1,...,d,
$$
$$
\bar{f}^{0}_{t}(x):=f^{0}_{t}(x+x_{t}),\quad
\bar{g}^{k}_{t}(x)=g^{k}_{t}(x+x_{t}),
$$
and the process
$x_{t}=(x^{1}_{t},...,x^{d}_{t})$ is defined by  
$$
x^{i}_{t}=-\int_{0}^{t}\sigma^{ik}_{s}(0)\,dw^{k}_{s}.
$$
This fact shows that the assertion of the present lemma
is a direct consequence of Lemma \ref{lemma 4.23.1}
in case the latter is applicable to \eqref{7.5.1}.

As is easy to see we will be able to apply Lemma \ref{lemma 4.23.1}
to \eqref{7.5.1} if we can find $\varepsilon'=
\varepsilon'(d,\delta,\varepsilon,p)>0$ such that
\begin{equation}
                                          \label{7.6.2}
\frac{1}{t-s}\int_{s}^{t}
(|\bar{a}^{ij}_{r}-\bar{a}^{ij}_{r(B)}|)_{(B)}\,dr
\leq\beta_{0},
\end{equation}
whenever $(s,t)\times B\in\bQ$ and $\rho(B)\leq\varepsilon'$.

 Denote by $N$, with or without subscripts,
various (large) constants depending only on $d $, $\delta$,
and $\varepsilon$
and observe that $|D\xi|\leq N $.
It follows easily that for
$B\in\bB$   we have
$$
(|\bar{a} ^{ij}_{r}-\bar{a} ^{ij}_{r(B )}|)_{(B )}
\leq(|\xi a^{ij}_{r}-(\xi a^{ij}_{r})_{(B+x_{r})}|)_{(B+x_{r})}
+\delta^{-1}\delta^{ij}
(|\xi  - \xi _{(B+x_{r})}|)_{(B+x_{r})}
$$
\begin{equation}
                                          \label{7.6.4}
\leq(|\xi a^{ij}_{r}-(\xi a^{ij}_{r})_{(B+x_{r})}|)_{(B+x_{r})}
+N_{1} \rho=:I_{r}+N_{1} \rho,
\end{equation}
where and below $\rho=\rho(B)$.

Let $z$ be the center of $B$ and set 
$$
y_{r}=(z+x_{r})(\rho+\varepsilon/2)
|z+x_{r}|^{-1}
$$
if $|z+x_{r}| \geq\rho+\varepsilon/2$ and $y_{r}=z+x_{r}$
otherwise. Observe that  $y_{r}$ is continuous in $r$
and  
\begin{equation}
                                          \label{7.6.3}
|y_{r}|\leq \rho+\varepsilon/2.
\end{equation}

Next we claim that
\begin{equation}
                                          \label{7.6.5}
I_{r}  
\leq  2 (|  a^{ij}_{r}-  a^{ij}_{r (B_{\rho}+y_{r})}|)_{(B_{\rho}+y_{r})}
+N_{2} \rho.
\end{equation}
  If \eqref{7.6.5} is true,
then by combining it with \eqref{7.6.4} and using \eqref{7.6.3}
we find that the left-hand side of \eqref{7.6.2} is less than
$$
(N_{1}+N_{2})\rho+2\sup_{|y|\leq \rho+\varepsilon/2}
\text{osc\,}(a^{ij},Q_{\rho}+(0,y),
0)
$$
if $\sigma_{t}^{nm}(0)=0$ for  all $t,n,m$ or, in general, less than
$$
(N_{1}+N_{2})\rho+2\text{Osc\,}(a^{ij},Q_{\rho},
\rho+\varepsilon/2),
$$
where $Q_{\rho}=(s,t)\times B_{\rho}$. Now
 \eqref{8.7.1} and \eqref{8.6.1}  imply 
that 
\eqref{7.6.2} is satisfied for $\rho\leq\varepsilon'$ if we choose
$\varepsilon'>0$ so that
$$
(N_{1}+N_{2})\varepsilon'\leq\beta_{0}/3,\quad
\varepsilon'\leq\varepsilon/4.
$$

Therefore, it only remains to prove the claim.
Obviously, if $|z+x_{r}|\geq\rho + \varepsilon/2$,
then $I_{r}=0$ and \eqref{7.6.5} holds.

In case $|z+x_{r}|<\rho +\varepsilon/2$   the estimates
$$
(|h_{r} -h_{ r(B')}|)_{(B')}\leq\frac{1}{|B'|^{2}}
\int_{B'}\int_{B'}|h_{r}(y)-h_{r}(z)|\,dy dz
\leq2(|h_{r} -h_{ r(B')}|)_{(B')},
$$
$$
|\xi(y)a^{ij}_{r}(y )-\xi(z)a^{ij}_{r}(z)|
\leq\xi(y )|a^{ij}_{r}(y )-a^{ij}_{r}(z)|
+N |\xi(y )-\xi(z)|
$$
show that
$$
I_{r}\leq
2(|  a^{ij}_{r}-  a^{ij}_{r (B+x_{r})}|)_{(B+x_{r})}
+N \rho ,
$$
which is equivalent to \eqref{7.6.5}. 
This proves the lemma in the particular case under consideration.

{\em General case\/}. We   rewrite the term 
$\Lambda^{k}_{t}u_{t}+g^{k}_{t}$ in \eqref{11.13.1}
as $\sigma^{ik}_{t}(0)\xi D_{i}u_{t}+\bar{g}^{k}_{t}$ with
$\bar{g}^{k}_{t}=g^{k}_{t}+(\sigma^{ik}_{t}-
\sigma^{ik}_{t}(0))D_{i}u_{t}$ and use the above result
to conclude that estimate 
\eqref{11.18.3} holds with $N=N_{1}=N_{1}(d,p,\delta,\varepsilon)$
if we add to its right-hand side
$$
N_{2}(d,p,\delta,\varepsilon)\beta_{1}\|Du\|_{\bL_{p}(\tau)}.
$$
By choosing $\beta_{1}=\beta_{1}(d,p,\delta,\varepsilon)$
so that $N_{2}\beta_{1}\leq1/2$, we get 
\eqref{11.18.3} with $2N_{1}$ in place of $N_{1}$.
The lemma is proved.

\begin{remark}   If
Assumptions \ref{assumption 1.2.1}
is satisfied with $K=0$
and   $a^{ij}_{t}$ and $\sigma^{ik}_{t}$
depend only on $\omega$ and $t$, then
the assertion of Lemma \ref{lemma 4.29.1} is true
with $\lambda_{0}=0$ and $N=N(d,p,\delta)$ and without requiring
$u$ to have compact support. This fact
can be obtained by following the arguments in Section 4.3
of \cite{Kr99}. Even though those arguments
are much longer, they allow one to prove a very general
result saying roughly speaking that ``whatever estimate
can be established
for solutions of the heat equation in Banach function spaces with norms
that are  
invariant under time dependent shifting of the $x$ coordinate,
the same estimate with the same constant also holds for solutions
of the parabolic equations with no lower order terms
and with the matrix of the second order coefficients
depending only on $t$ and dominating (in the matrix sense)
the unit matrix"
(see \cite{Kr94_1}).

\end{remark}

Next step is to consider equations with
lower order terms.
The following lemma and its corollary are stated
in a slightly more general form than it is needed in the present
article. The point is that we intend to use them
in a subsequent article about equations in half spaces.

\begin{lemma}
                                    \label{lemma 12.7.1}
Let $G\subset\bR^{d}$ be a domain (perhaps, $G=\bR^{d}$)
  and take   $\varepsilon\geq
\varepsilon_{1}>0$ and
 $\varepsilon_{2}\in(0,\varepsilon/4]$.

(i) Let  
$f^{j},g\in\bL_{p}(\tau)$ and
let $u\in\cW^{1}_{p,0}(\tau)$
satisfy~\eqref{11.13.1} for $t\leq\tau$
and be such that 
$u_{t}(x)=0$ if $x\not\in G$.

(ii)  Suppose that Assumptions \ref{assumption 1.2.1}
is satisfied.  

(iii) Suppose that assumption (iii) of Lemma
\ref{lemma 4.29.1} is satisfied for any $t_{0}\geq0$
and $x_{0} $ such that $\text{\rm dist}\,(x_{0},G)\leq\varepsilon_{2}$.

Then there exist constants  $N, \lambda_{0}\geq0$, 
depending only on $d$, $p$, $K$, $\delta$, $\varepsilon$, 
$\varepsilon_{1}$, and  
$\varepsilon_{2}$,
such that estimate \eqref{11.18.3} holds true
whenever $\lambda\geq\lambda_{0}$.
 
\end{lemma}

Proof.  
As usual we will use partitions of unity.
Take a nonnegative
$\xi\in C^{\infty}_{0}(B_{\varepsilon_{2}})$
with unit $L_{p}$-norm and  
take a nonnegative
$\eta\in C^{\infty}_{0}((0,\varepsilon_{1}^{2}))$
with unit $L_{p}$-norm.
For $s\in\bR$ and $y\in\bR^{d}$
introduce
$$
\zeta(t,x)=\xi(x)\eta(t),\quad
\zeta^{s,y}(t,x)=\zeta(t-s,x-y),\quad u^{s,y}_{t}(x)
=\zeta^{s,y}(t,x)u_{t}(x)
$$ 
so that, in particular,
\begin{equation}
                                           \label{12.7.1}
|u_{t}(x)|^{p}=\int_{\bR^{d+1}}|u^{s,y}_{t}(x)|^{p}\,dyds.
\end{equation}

Observe that for each $s,y$
$$
du^{s,y}_{t}=\big(\sigma^{ik}_{t} D_{i}u^{s,y}_{t}
+\hat{g}^{s,y,k}_{t}\big)\,dw^{k}_{t} 
$$
\begin{equation}
                                           \label{11.13.3}
+\big(D_{j}(a^{ij}_{t}D_{i}u^{s,y}_{t}) 
-\lambda u^{s,y}_{t}
+D_{j}\hat{f}^{s,y,j}_{t}+\hat{f}^{s,y,0}_{t}\big)\,dt
\end{equation}
for $t\leq\tau$,
where we dropped the argument $x$ (and $\omega$) and
$$
\hat{g}^{s,y,k}_{t}=\zeta^{s,y}(\nu^{k}_{t}u_{t}+g^{k}_{t})
-u_{t}\sigma^{ik}_{t}D_{i}\zeta^{s,y},
$$
$$
\hat{f}^{s,y,j}_{t}=\zeta^{s,y}(a^{j}_{t}u_{t}+f^{j}_{t})
-a^{ij}_{t}u_{t}D_{i}
\zeta^{s,y} ,\quad j=1,...,d,
$$
$$
\hat{f}^{s,y,0}_{t}=\zeta^{s,y}(f^{0}_{t}+b^{i}_{t}D_{i}u_{t}
+c_{t}u_{t})-f^{j}_{t}
D_{j}\zeta^{s,y}-(a^{ij}_{t}D_{i}u_{t}+a^{j}_{t}u_{t})D_{j}\zeta^{s,y}
+\zeta^{s,y}_{t}u_{t},
$$
and $\zeta^{s,y}_{t}(t,x)=\xi(x-y)\eta'(t-s)$.

As is easy to see
 $u^{s,y}(t,x)=0$ for $(t,x)\not\in 
(s_{+},s_{+}+\varepsilon_{1}^{2})\times B_{\varepsilon_{2}}
(y)$. Therefore, by Lemma \ref{lemma 4.29.1} 
if $\text{\rm dist}\,(y,G)\leq\varepsilon_{2}$, then
$$
 \lambda^{p/2}\|u^{s,y}\|^{p}_{\bL_{p}(\tau)} 
+\|Du^{s,y} \|^{p}_{\bL_{p}(\tau)} 
$$
\begin{equation}
                                             \label{8.7.5}
\leq 
 N\big(\sum_{j=1}^{d}
\|\hat{f}^{s,y,j}\|^{p}_{\bL_{p}(\tau)}
+\|\hat{g}^{s,y}\|^{p}_{\bL_{p}(\tau)}
\big) +N\lambda^{-p/2}\|\hat{f}^{s,y,0}\|^{p}_{\bL_{p}(\tau)}
\end{equation}
provided that $\lambda\geq\lambda_{0}$, where $N$ and $\lambda_{0}$
depend only on $d,\delta,p$, and $\varepsilon$.
This estimate also, obviously, holds if
$\text{\rm dist}\,(y,G)>\varepsilon_{2}$ since then
$u^{s,y}_{t}\equiv0$.

Next,
$$
|\hat{f}^{s,y,j}_{t}|\leq  
N \bar{\zeta}^{s,y} |u_{t}| +\zeta^{s,y}|f^{j}_{t}|,\quad j=1,...,d,
$$
$$
|\hat{f}^{s,y,0}_{t}|\leq N \bar{\zeta}^{s,y}( |Du_{t}|
+|u_{t}|)+N \bar{\zeta}^{s,y}\sum_{j=0}^{d}|f^{j}_{t}|,
$$
$$
|\hat{g}^{s,y}_{t}|_{\ell_{2}}\leq 
 N\bar{\zeta}^{s,y}|u_{t}|+\zeta^{s,y}|g_{t}|_{\ell_{2}},
$$
where  
$\bar{\zeta}=\zeta+|D\zeta|+
|\zeta_{t}|$, $
\bar{\zeta}^{s,y}(t,x)=\bar{\zeta}(t-s,x-y)$, and here and below
we allow
the constants $N $ to depend only on $d,p,\delta$,  $K$, $\varepsilon$,
  $\varepsilon_{1}$, and $\varepsilon_{2}$. 

We also  notice that $|\zeta^{s,y}Du_{t}|\leq
|D(\zeta^{s,y}u_{t})|+\bar{\zeta}^{s,y}|u_{t}|$. Then we find that
$$
\lambda^{p/2}\|\zeta^{s,y}u\|^{p}_{\bL_{p}(\tau)}
+\| \zeta^{s,y}Du \|^{p}_{\bL_{p}(\tau)}
$$
$$
\leq  N\big(\sum_{i=1}^{d}
\|\bar{\zeta}^{s,y}f^{i}\|^{p}_{\bL_{p}(\tau)}+\|\zeta^{s,y}
g\|^{p}_{\bL_{p}(\tau)}+\|\bar{\zeta}^{s,y}
u\|^{p}_{\bL_{p}(\tau)}
\big)
$$
$$
+N\lambda^{-p/2}(
\|\bar{\zeta}^{s,y}f^{0}\|^{p}_{\bL_{p}(\tau)}
+\|\bar{\zeta}^{s,y}Du\|^{p}_{\bL_{p}(\tau)}).
$$
We integrate through this estimate and use
formulas like \eqref{12.7.1}. Then we obtain 
$$
\lambda^{p/2}\| u\|^{p}_{\bL_{p}(\tau)}
+\|D u \|^{p}_{\bL_{p}(\tau)}
$$
$$
\leq  N_{1} \big(\sum_{i=1}^{d}
\| f^{i}\|^{p}_{\bL_{p}(\tau)}+\| 
g\|^{p}_{\bL_{p}(\tau)}+\| u\|^{p}_{\bL_{p}(\tau)}
\big)+N_{1} \lambda^{-p/2}(\| f^{0}\|^{p}_{\bL_{p}(\tau)}
+\| Du\|^{p}_{\bL_{p}(\tau)}).
$$
Finally, we  increase $\lambda_{0}\geq0$,
if necessary,
in such a way that
$N_{1}\lambda^{-p/2}\leq1/2$
for $\lambda\geq\lambda_{0}$.
 Then we obviously 
arrive at \eqref{11.18.3} with
$N=2N_{1}$. The lemma is proved.

To the best of the author's knowledge
the following multiplicative estimate
is new even in the deterministic case.

\begin{corollary}
                                    \label{corollary 12.11.1}
Let $\lambda=0$. Then under the assumptions of Lemma 
\ref{lemma 12.7.1} we have  
$$
 \|Du \|_{\bL_{p}(\tau)} 
\leq N\big(\sum_{i=1}^{d}
\|f^{i}\|_{\bL_{p}(\tau)}+\|g\|_{\bL_{p}(\tau)}
+\|f^{0} \|^{1/2}_{\bL_{p}(\tau)}\|u\|_{\bL_{p}(\tau)}^{1/2}
+\|u\|_{\bL_{p}(\tau)}
\big),
$$
where $N$ depends only on $d,p,K,\delta$, 
 $\varepsilon$, 
$\varepsilon_{1}$, and  
$\varepsilon_{2}$.

\end{corollary}

Indeed,    take a $\lambda>0$ and
 add and subtract the term 
$(\lambda_{0}+\lambda)u_{t}\,dt$ on the right in
\eqref{11.13.1}, thus introducing
$\lambda$ into the equation and modifying $f^{0}_{t}$
by including into it one of $(\lambda_{0}+\lambda)u_{t}$.
Then after  applying \eqref{11.18.3},
 we see that
$$
 \|Du \|_{\bL_{p}(\tau)} 
\leq N\big(\sum_{i=1}^{d}
\|f^{i}\|_{\bL_{p}(\tau)}+\|g\|_{\bL_{p}(\tau)}
$$
$$
+(\lambda_{0}+\lambda)^{-1/2}
\|f^{0} \| _{\bL_{p}(\tau)} 
+(\lambda_{0}+\lambda)^{1/2}\|u\|_{\bL_{p}(\tau)}
\big).
$$
Now it only remains to take the inf with respect to
$\lambda>0$. 

{\bf Proof of Lemma \ref{lemma 2.25.1}}. By 
bearing in mind an obvious shifting of  time
we see that in the proof of assertions (i)-(iii)
we may assume that $s=0$.

(i) First of all observe that uniqueness of
solution of \eqref{2.27.2} is well known
even in a much wider class than $\cW^{1}_{p}(\infty)$.

Let $u_{0}\in\tr_{\!0}\cW^{1}_{p}$, then 
$u_{0}\in W^{1-2/p}$ for almost each $\omega$  
and there is a unique solution
of the heat equation
$$
dv_{t}=\Delta v_{t}\,dt
$$
of class $L_{p}((0,1),W^{1}_{p})$ with initial 
condition $u_{0}$. 
Furthermore,
$$
\|v\|_{L_{p}((0,1),W^{1}_{p})}
\sim\|u_{0}\|_{W^{1-2/p}_{p}}.
$$

Next take a $\zeta\in C^{\infty}_{0}(\bR)$
such that $\zeta_{0}=1$ and 
$\zeta_{t} =0$ for $t\geq1/2$
and define $\psi_{t}(x)=e^{-t}v_{t}(x)\zeta_{t}$
 for $t\in[0,1]$
and as zero if $t\geq1/2$. Notice that
 (a.s.)
$$
\psi\in L_{p}(\bR _{+},W^{1}_{p}),
$$ 
and
$$
\frac{\partial}{\partial t}\psi_{t}
=\Delta\psi_{t}-\psi_{t}+e^{-t}\zeta'_{t}v_{t}
$$

Then it is a classical result
that there exists a unique $\phi
\in L_{p}(\bR _{+},W^{2}_{p}) $
which solves the equation
$$
d\phi_{t}=(\Delta \phi_{t}-\phi_{t} +e^{-t}\zeta'_{t}v_{t})\,dt
$$
with {\em zero\/} initial condition. In addition,
$$
\|\phi\|_{L_{p}(\bR _{+},W^{2}_{p})}\leq
 N\|\zeta'v\|_{L_{p}(\bR _{+},L_{p})}   
\leq N\|u_{0}\|_{W^{1-2/p}_{p}},
$$
where the constants $N$ depend only on $d$ and $p$.
Owing to these estimates and uniqueness,
the operators mapping
$u_{0}$ into $v$ and $\phi$ are continuous
(and nonrandom).
Since $u_{0}$ is
$\cF_{0}$-measurable, the same is true for
$\psi$, $\phi$, and $u=\psi-\phi$,
which is of class $L_{p}((0,1),W^{1}_{p})$,
satisfies \eqref{2.27.2}
and equals $u_{0}$ for $t=0$. Also
 for each $\omega$
$$
\|u\|_{L_{p}(\bR_{+},W^{1}_{p})}\leq
\|\psi\|_{L_{p}(\bR_{+},W^{1}_{p})}+
\|\phi\|_{L_{p}(\bR_{+},W^{1}_{p})}\leq
N\|u_{0}\|_{W^{1-2/p}_{p}},
$$
where $N$ depends only on $d$ and $p$. By 
raising the extreme terms to the $p$th power and taking
expectations we get the first inequality 
in \eqref{2.25.4} and also finish proving the ``only
if" part of (i).  
 
To prove the ``if" part assume that we have a $v
\in\cW^{1}_{p}(\infty)$ satisfying \eqref{2.27.2}  
and equal $u_{0}$ at $t=0$. Then
 $u_{t}=v_{t}e^{t}$ satisfies
$\partial u_{t}/\partial t=\Delta u_{t}$ and is of
class $\cW^{1}_{p}(1)$. It follows that 
almost all $\omega$ we have $u\in 
L_{p}((0,1),W^{1}_{p})$,  $u_{0}\in W^{1-2/p}_{p}$,
and
$$
\|u_{0}\|_{W^{1-2/p}_{p}}\leq N\|u\|_{L_{p}((0,1),W^{1}_{p})}
\leq N\|v\|_{L_{p}(\bR_{+},W^{1}_{p})}.
$$
By raising all expressions to the power $p$
and taking expectations we arrive at the second
estimate in \eqref{2.25.4}. Assertion (i) is proved. 

The ``only if" part in (ii) is, actually, proved above.
To prove the ``if" part  write
$$
dv_{t}=(D_{i}f^{i}_{t}+f^{0}_{t})\,dt
+g^{k}_{t}\,dw^{k}_{t}=(\Delta v_{t}-\lambda v_{t}
+D_{i}\hat{f}^{i}_{t}+\hat{f}^{0}_{t})\,dt
+g^{k}_{t}\,dw^{k}_{t},
$$
where the constant $\lambda>0$
will be chosen later,
$\hat{f}^{i}_{t}=f^{i}_{t}-D_{i}v_{t}$, $i=1,...,d$,
$\hat{f}^{0}_{t}=f^{0}_{t}+\lambda v_{t}$, and
$\hat{f}^{j},g\in \bL_{p}(1)$. 
Next, take the function $\zeta$
as above, set $u=v\zeta$, and observe that  
\begin{equation}
                                               \label{2.29.1}
du_{t}=(\Delta u_{t}-\lambda u_{t}
+D_{i} \check{f}^{i}_{t}+\check{f}^{0}_{t})\,dt
+\check{g}^{k}_{t}\,dw^{k}_{t},
\end{equation}
where $\check{f}^{0}=\zeta\hat{f}^{0}+v\zeta'$,
$\check{f}^{i}_{t}=\zeta \hat{f}^{i}_{t}$, $i=1,...,d$,
$\check{g}^{k}=\zeta g^{k}$
and $\check{f}^{j},\check{g}\in \bL_{p}(\infty)$
and $u\in\cW^{1}_{p}(\infty)$.

By Lemma \ref{lemma 4.23.1}, for $\lambda$ fixed and large enough
(actually, one can take $\lambda=1$, which is shown by using
dilations),
equation \eqref{2.29.1} with zero initial condition
admits a unique solution $\psi\in\cW^{1}_{p}(\infty)$
and 
$$
\|\psi\|_{\bW^{1}_{p}(\infty)}
\leq N(\sum_{j=0}^{d}
\|\check{f}^{j}\|_{\bL_{p}(\infty)}
+\|\check{g}\|_{\bL_{p}(\infty)})
$$
$$
\leq N(\sum_{j=0}^{d}
\|f^{j}\|_{\bL_{p}(1)}
+\|g\|_{\bL_{p}(1)}+\|v\|_{\bW^{1}_{p}(1)}).
$$
Then the difference $\phi=u-\psi$ satisfies
\eqref{2.27.2}, is of class $\cW^{1}_{p}(\infty)$,  
and $\phi_{0}=u_{0}$. By assertion (i)
we have $u_{0}\in\tr_{\!0}\cW^{1}_{p}$,
which proves the ``if" part in (ii). Furthermore,
$$
\|u_{0}\|_{\tr_{\!0}\cW^{1}_{p}} \leq 
N \|\phi\|_{\bW^{1}_{p}(\infty)}\leq
N \|u\|_{\bW^{1}_{p}(\infty)}+
N \|\psi\|_{\bW^{1}_{p}(\infty)}
$$
$$
\leq
N \|v\|_{\bW^{1}_{p}(1)}+
N \|\psi\|_{\bW^{1}_{p}(\infty)}
\leq N(\sum_{j=0}^{d}
\|f^{j}\|_{\bL_{p}(1)}
+\|g\|_{\bL_{p}(1)}+\|v\|_{\bW^{1}_{p}(1)}).
$$
This   proves 
assertion (iii). 

To prove (iv) observe that obvious dilations of the $t$ axis
allow us to assume that $s=1$. Then
 write \eqref{12.3.1} for $t\in[0,1]$ and
notice that $tu_{t}$ admits representation
\eqref{12.3.1} with new $f^{j}$ and $g^{k}$
having simple relations with $u_{t}$ and the original
$f^{j}$ and $g^{k}$. It follows that in the rest of the 
proof we may assume that $u_{0}=0$.

In that case take a sufficiently large $\lambda>0$
and consider the equation
$$
d v_{t}=(\Delta v_{t}-\lambda v_{t}+D_{i}\bar{f}^{i}_{t}
+\bar{f}^{0}_{t})\,dt+\bar{g}^{k}_{t}\,dw^{k}_{t}
$$
for $t\geq0$ with zero initial condition,
where
$$
\bar{f}^{i}_{t}=f^{i}_{t}I_{(0,1)}(t)
-D_{i}u_{t}I_{(0,1)}(t),\quad i=1,...,d,
$$
$$
\bar{f}^{0}_{t}=( f^{0}_{t}+\lambda u_{t})I_{(0,1)}(t),\quad
\bar{g}^{k}_{t}=g^{k}_{t}I_{(0,1)}(t).
$$
By uniqueness, $ v_{t}=u_{t}$ for $t\in[0,1]$ and by assertion (iii)
we have $ v_{1}\in\tr_{\!1}\cW^{1}_{p}$. This fact combined
with already known estimates of $v$ proves assertion (iv).
The lemma is proved.

\mysection{Proof of Theorem \protect\ref{theorem 12.7.1}}
                                    \label{setion 3.2.1}

Owing to Lemma \ref{lemma 2.25.1} we may assume that
we are given a $v$ as in assertion (i) of the lemma.
By introducing a new unknown function $\bar{u}=u-v$
we   see that $u$ satisfies 
\eqref{11.13.1} and $u_{0}=v_{0}$ if and only if
$\bar{u}_{0}=0$ and
$$
d\bar{u}_{t}=(L_{t}\bar{u}_{t}-\lambda\bar{u}_{t}
+D_{j}\bar{f}^{j}_{t}+\bar{f}^{0}_{t})\,dt+
(\Lambda^{k}_{t}\bar{u}_{t}+\bar{g}^{k}_{t})\,dw^{k}_{t},
$$
where
$$
\bar{f}^{j}_{t}=f^{j}_{t}-D_{j}v_{t} 
+a^{ij}_{t}D_{i}v_{t}+a^{j}_{t}v_{t},\quad j=1,...,d,
$$
$$
\bar{f}^{0}_{t}=f^{0}_{t} 
+b^{i}_{t}D_{i}v_{t}+(c_{t}-\lambda+1)v_{t},
$$
$$
\bar{g}^{k}_{t}=g^{k}_{t} 
+\sigma^{ik}_{t}D_{i}v_{t}+\nu^{k}_{t}v_{t}.
$$
By Lemma \ref{lemma 2.25.1} we have
 $\bar{f}^{j},\bar{g}\in\bL_{p}(\tau)$
and the problem of finding solutions of \eqref{11.13.1}
with initial data $u_{0}$ is thus reduced to the same problem
but with zero initial data.

Furthermore, if estimate \eqref{1.4.2} holds for solutions
with zero initial condition, then (for $\lambda\geq\lambda_{0}$)
$$
 \lambda^{1/2}\|u\|_{\bL_{p}(\tau)} 
+\|Du \|_{\bL_{p}(\tau)} 
- \lambda^{1/2}\|v\|_{\bL_{p}(\tau)} 
-\|Dv \|_{\bL_{p}(\tau)} 
$$
$$
\leq
\lambda^{1/2}\|\bar u\|_{\bL_{p}(\tau)} 
+\|D\bar u \|_{\bL_{p}(\tau)} 
$$
$$
\leq N\big(\sum_{i=1}^{d}
\|\bar f^{i}\|_{\bL_{p}(\tau)}+\|\bar g\|_{\bL_{p}(\tau)}
\big)+N\lambda^{-1/2}\|\bar f^{0}\|_{\bL_{p}(\tau)}
$$
$$
  \leq N\big(\sum_{i=1}^{d}
\|f^{i}\|_{\bL_{p}(\tau)}+\|g\|_{\bL_{p}(\tau)}
+\| v\|_{\bW^{1}_{p}(\tau)}
\big)
$$
$$
+N\lambda^{-1/2}(\|f^{0}\|_{\bL_{p}(\tau)}+
 \| v\|_{\bW^{1}_{p}(\tau)})+N
\lambda^{1/2}\| v\|_{\bL_{p}(\tau)},
$$
which yields \eqref{1.4.2} in full generality.

It follows that while proving \eqref{1.4.2} we may
also assume that $u_{0}=0$. Therefore, in the rest of the 
proof of assertions (i) and (ii) we assume that $u_{0}=0$.
 Having in mind the substitution $u_{t}=v_{t}e^{-\mu t}$,
we see that while proving assertion (i) it suffices to concentrate
on large $\lambda$ and prove only the second part of the assertion.

We recall that we suppose that
Assumption \ref{assumption 1.2.6} 
is satisfied with $\beta=\beta_{0}/3$
and $\beta_{0}$ from Lemma \ref{lemma 4.23.1} 
and Assumption \ref{assumption 8.9.1} 
is satisfied with
$\beta_{1}$
defined in Lemma \ref{lemma 4.29.1}. It follows that
assumption (iii) of Lemma \ref{lemma 4.29.1} 
is satisfied for any $(t_{0},x_{0})$.

Now we take $\lambda_{0}$
larger than the one in Lemma \ref{lemma 12.4.1}
and the one in Lemma \ref{lemma 12.7.1}. In that case
uniqueness follows from Lemma \ref{lemma 12.4.1}.
In the proof of existence
we will rely on the method of continuity
and the a priori estimate \eqref{11.18.3} which
is established in Lemma \ref{lemma 12.7.1}.
For $\lambda\geq\lambda_{0}$ and $\theta\in[0,1]$
we consider the equation
\begin{equation}
                                          \label{12.7.3}
du_{t}=[(\theta L_{t}+(1-\theta)\Delta)u_{t}
-\lambda u_{t}+D_{i}f^{i}_{t}+f^{0}_{t})\,dt
+(\theta \Lambda^{k}_{t}u_{t}+g^{k}_{t})\,dw^{k}_{t}.
\end{equation}
 
 We call a $\theta\in[0,1]$ ``good" if
the assertions of the theorem hold for
equation \eqref{12.7.3}. Observe that $0$
is a ``good" point by Lemma \ref{lemma 4.23.1}.
Now to prove the theorem it suffices to show
that there exists a  $\gamma>0$ such that
if  $\theta_{0}$ is a good point then
all points of the interval $[\theta_{0}-\gamma,
\theta_{0}+\gamma]\cap[0,1]$ are ``good".
So fix a ``good" $\theta_{0}$ and 
for any $v\in\bW^{1}_{p}(\tau)$ consider the equation
$$
du_{t}=[(\theta_{0} L_{t}+(1-\theta_{0})\Delta)u_{t}
-\lambda u_{t}
+(\theta-\theta_{0})(L_{t}-\Delta)v_{t}
+D_{i}f^{i}_{t}+f^{0}_{t})\,dt
$$
\begin{equation}
                                          \label{12.7.4}
+(\theta_{0} \Lambda^{k}_{t}u_{t}
+(\theta-\theta_{0})\Lambda^{k}v_{t}+g^{k}_{t})\,dw^{k}_{t}.
\end{equation}
Observe that
$$
(L_{t}-\Delta)v_{t}=D_{j}\big((a^{ij}-\delta^{ij})
D_{i}v_{t}+a^{j}_{t}v_{t}\big)+b^{i}_{t}D_{i}v_{t}
+cv_{t}
$$
and recall that $v\in\bW^{1}_{p}(\tau)$. It follows  by assumption that
equation \eqref{12.7.4} has a unique
solution $u\in\cW^{1}_{p,0}(\tau)$ ($\subset
\bW^{1}_{p }(\tau)$). 

In this way, for $f^{j}$ and $g$ being fixed,
 we define a mapping
$v\to u$ in the space $\bW^{1}_{p}(\tau)$. 
It is important to keep in mind
 that the image $u$ of 
$v\in\bW^{1}_{p}(\tau)$ is always in 
$\cW^{1}_{p,0}(\tau)$.
Take
$v',v''\in\bW^{1}_{p}(\tau)$ and let $u',u''$ be their 
corresponding images. Then $u:=u'-u''$ satisfies
$$
du_{t}=[(\theta_{0} L_{t}+(1-\theta_{0})\Delta)u_{t}
-\lambda u_{t}
+(\theta-\theta_{0})(L_{t}-\Delta)v_{t})\,dt
$$
$$
+(\theta_{0} \Lambda^{k}_{t}u_{t}
+(\theta-\theta_{0})\Lambda^{k}v_{t} )\,dw^{k}_{t},
$$
where $v=v'-v''$. It follows by Lemma \ref{lemma 12.7.1}
that
$$
\|u\|_{\bW^{1}_{p}(\tau)}\leq N|\theta-\theta_{0}|
\,\|v\|_{\bW^{1}_{p}(\tau)}
$$
with a constant $N$ independent of $v'$, $v''$, $\theta_{0}$,
and $\theta$. 
For $\theta$ sufficiently close
to $\theta_{0}$, our mapping is a contraction
and, since $\bW^{1}_{p}(\tau)$ is a Banach
space, 
it  has a fixed point. This fixed point is in $\cW^{1}_{p,0}
(\tau)$ and,
obviously, satisfies \eqref{12.7.3}. This proves
assertion  (i)   of the theorem.

Estimate \eqref{1.4.2} is proved above in Lemma \ref{lemma 12.7.1}
and assertion (iii) follows from Theorem \ref{theorem 12.3.1}.
The theorem is proved.

\end{document}